\documentclass[a4paper,11pt]{elsarticle}
\usepackage{amsmath,amssymb,amsthm,graphicx,subfigure,float,caption,epstopdf,tabularx,color, bm,amsfonts,epic}
\usepackage[top=1in, bottom=1in, left=1.25in, right=1.25in]{geometry}
\usepackage{appendix}
\usepackage{multirow}
\usepackage{diagbox}
\newtheorem{thm}{Theorem}[section]

\newtheorem{rmk}{Remark}[section]
\newtheorem{shm}{Scheme}[section]
\newtheorem{ex}{Example}[section]
\newtheorem*{prf}{Proof}
\numberwithin{equation}{section}

\usepackage{indentfirst}
\graphicspath{{Fig}}

\usepackage{hyperref}
\hypersetup{
	colorlinks=true,
	citecolor=blue,
	linkcolor=red,
	urlcolor=green}

\newcommand{\be}{\begin{equation}}
\newcommand{\ee}{\end{equation}}
\newcommand{\ba}{\begin{array}}
\newcommand{\ea}{\end{array}}
\newcommand{\bea}{\begin{eqnarray}}
\newcommand{\eea}{\end{eqnarray}}
\newcommand{\beas}{\begin{eqnarray*}}
\newcommand{\eeas}{\end{eqnarray*}}

\begin{document}

\begin{frontmatter}

\title{High-order structure-preserving schemes for the regularized logarithmic Schr\"{o}dinger equation}
\author{Fan Yang}
\ead{imqingfeng17@163.com}

\author{Zhida Zhou}
\ead{1415747150@qq.com}

\author{Chaolong Jiang\corref{4}}
\ead{chaolong\_jiang@126.com}
\address{School of Statistics and Mathematics,Yunnan University of Finance and Economics, Kunming 650221, China}

\cortext[4]{Corresponding author.}

\begin{abstract}
In this paper, a novel high-order, mass and energy-conserving scheme is proposed for the regularized logarithmic Schr\"{o}dinger equation(RLogSE). Based on the idea of the supplementary variable method (SVM), we firstly reformulate the original system into an equivalent form by introducing two supplementary variables, and the resulting SVM reformulation is then discretized by applying a high-order prediction-correction scheme in time and a Fourier pseudo-spectral method in space, respectively. The newly developed scheme can produce numerical solutions along which the mass and original energy are precisely conserved, as is the case with the analytical solution. Additionally, it is extremely efficient in the sense that only requires solving a constant-coefficient linear systems plus two algebraic equations, which can be efficiently solved by the Newton iteration at every time step. Numerical experiments are presented to confirm the accuracy and structure-preserving properties of the new scheme.  \\[2ex]
\end{abstract}

\end{frontmatter}

\tableofcontents

\section{Introduction}
The logarithmic Schr\"{o}dinger equation (LogSE) has been widely used in different branches of fundamental physics, such as nuclear physics~\cite{Hefter1985}, diffusion phenomena~\cite{Hansson2009} and Bose-Einstein condensation~\cite{Zloshchastiev2019}. In this paper, we consider the following logarithmic Schr\"{o}dinger equation
\begin{align}\label {LogSE}
\left\{
\begin{matrix}
\mathrm{i}u_{t}(\textbf{x}, t) + \Delta u(\textbf{x}, t) = \lambda u(\textbf{x}, t)\ln(|u(\textbf{x}, t)|^{2}), \ t > 0,\ \textbf{x} \in \Omega\subset\mathbb{R}^d,\ d=1,2,  \hfill\hfill\\
u(\textbf{x}, 0) = u_{0}(\textbf{x}), \textbf{x} \in \bar\Omega\subset\mathbb{R}^d, \hfill\hfill\\
\end{matrix}
\right.
\end{align}
where~$u(\textbf{x},t)$ is the wave function,~\textbf{x}~is the spatial variable, ~$t$~is time variable,~$\lambda\in \mathbb{R}/\{0\}$ is a dimensionless real constant of the nonlinear interaction strength,~$\Omega$ is a bounded domain with a periodic boundary condition.
The LogSE \eqref{LogSE} conserves the following two invariants~\cite{2}:
\ the mass
\begin{align}
\mathcal{M}[u(\cdot,t)]= \int_{\Omega} |u|^{2} \, d\textbf{x} \equiv \mathcal{M}[u(\cdot,0)],
\end{align}
and energy
\begin{align}
\mathcal{E}[u(\cdot,t)] := \int_{\Omega} (|\nabla u|^{2} + \lambda |u|^{2} \ln (|u|^{2})) \, d\textbf{x} \equiv \mathcal{E}[u(\cdot,0)].
\end{align}

To overcome the singularity of the logarithmic nonlinearity at the origin, a regularized logarithmic Schr\"{o}dinger equation (RLogSE) model with a small regularized parameter $0<\varepsilon\ll1$ was introduced \cite{Bao2019}, as follows:
\begin{align}\label {RLogSE}
\left\{
\begin{matrix}
\mathrm{i}u_{t}^{\varepsilon}(\textbf{x},t)+\Delta u^{\varepsilon}(\textbf{x},t)=\lambda u^{\varepsilon}(\textbf{x},t)\ln(\varepsilon +|u^{\varepsilon}(\textbf{x},t)|)^{2},  t > 0, \textbf{x} \in \Omega\subset\mathbb{R}^d,\ d=1,2, \hfill\hfill\\
u^{\varepsilon}(\textbf{x},0)=u_{0}(\textbf{x}),\textbf{x}\in \bar{\Omega}\subset\mathbb{R}^d.\hfill\hfill\\
\end{matrix}
\right.
\end{align}
It is shown that the RLogSE \eqref{RLogSE} approximates the LogSE \eqref{LogSE} with linear convergence rate $ \mathcal{O}(\epsilon)$~\cite{BCST-NM-2019},
and has mass and energy conservation laws~\cite{Bao2019} similar to those of the original model, i.e., the mass
\begin{align}
\label {regularized mass}
\mathcal{M}^{\epsilon}[u^{\epsilon}(\cdot,t)] := \int_{\Omega} |u^{\epsilon}|^{2} \, d\textbf{x} \equiv \mathcal{M}^{\epsilon}[u^{\epsilon}(\cdot,0)], 
\end{align}
the regularized energy
\begin{align}
\label {regularized energy}
\mathcal{E}^{\epsilon}[u^{\epsilon}(\cdot,t)] := \int_{\Omega} \left[ |\nabla u^{\epsilon}|^{2} + 2\epsilon \lambda |u^{\epsilon}| + 2\lambda (|u^{\epsilon}|^{2} - \epsilon^{2}) \ln(\epsilon + |u^{\epsilon}|)
\right ] d\textbf{x}\equiv \mathcal{E}^{\epsilon}[u^{\epsilon}(\cdot,0)].
\end{align}

In the last few decades, there has been considerable literature on theoretical and numerical studies on the classical  nonlinear Schr\"{o}dinger equation. For more details, please refer to the review papers \cite{ABC-CPC-2013,BBGI-JCAM-2024} and the references therein. Nevertheless, since the blow up of the logarithmic nonlinearity, there are significant difficulties in designing numerical methods for LogSE \eqref{LogSE}. Thus, the exploration of nonlinear Schr\"{o}dinger equations incorporating logarithmic terms remains nascent. Scott and Shertzer~\cite{Scott2018} proposed an iterative finite element method to solve the Coulombic logSE for the spherically symmetric states. Later on, Bao and his collaborators \cite{Bao2019} first introduced the idea of the regularized LogSE model \eqref{RLogSE} with a regularization parameter. Then, based on the regularized LogSE model \eqref{RLogSE}, they present and analyze different regularized numerical methods \cite{Bao2019,BCST-NM-2019,BCST-MMMAS-2022}. Other works on well-posedness and numerical methods of the regularized LogSE model can be found in Refs. \cite{CS-arXiv-2022,CG-ANM-2021,LZH-ANM-2019}. Different from existing works based on regularisation, Paraschis and Zouraris proposed and analyze an implicit Crank-Nicolson finite difference scheme without regularizing the logarithmic term~\cite{Paraschis2023}. Then, Wang et al.\cite{LYZ2024} presented and analysed a first-order IMEX scheme of the LogSE without regularization. More recently, a direct (non-regularised) Lie-Totter time-splitting Fourier spectral scheme with low regularity initial data and solution is analyzed in \cite{ZW-arXiv-2024}. Among these numerical methods, most of existing works mainly focused on the stability and convergence of the numerical schemes, while ignoring the mass and energy conservation laws of the model, which is crucial for maintaining the stability and accuracy of the simulation in long computations.

In \cite{BCST-NM-2019,Paraschis2023}, the mass- or energy-conserving Crank-Nicolson schemes are presented and analyzed. However, they are fully implicit, so that at every time step, a nonlinear equation shall be solved by using a nonlinear iterative method and thus it may be time consuming. Additionally, these schemes are second-order accurate in time, which usually  produces significantly large numerical error for a given large time step. More recently, based on the energy quadratization (EQ) approach \cite{SXY-JCP-2018,YZW-JCP-2017}, Qian et al.~\cite{qian2023} developed a class of regularized high-order numerical schemes for the RLogSE \eqref{RLogSE} that preserve the mass and energy conservation laws. Despite the proposed scheme shows excellent numerical behaviors in long computations, it is fully implicit and only conserves a modified energy rather than the original one.

Due to the high computational cost of the fully implicit schemes in the literatures, ones are devoted to construct linearly implicit energy-conserving schemes, in which a linear system is to be solved at every time step. Thus it is computationally much cheaper than that of the implicit scheme. In \cite{DO-SISC-2011,FM2011}, Furihata, Owren et al. presented a general framework for deriving linearly implicit energy-conserving schemes of partial difference equations with polynomial nonlinearities. In \cite{CLW18}, Cai et al. presented the partitioned averaged vector field method, which provides an efficient approach to design linearly implicit mass- and energy-preserving for the Klein-Gordon-Schr\"odinger equation. However, numerical schemes bases on these approach are second-order accurate in time at most, and invalid for the LogSE \eqref{RLogSE}. Actually, based on the theories on the EQ approach and quadratic invariants-conserving Rung-Kutta (RK) methods \cite{HLWbook2006}, we can easily proposed high-order, linearly-implicit energy-conserving schemes of the LogSE \eqref{RLogSE}. Nevertheless, as pointed out above, the proposed schemes only conserve a modified energy. We note that the time-splitting scheme is mass-conserving, unconditionally stable and easy to be high order, which however they cannot conserve the Hamiltonian energy \cite{HLWbook2006,CG-ANM-2021}. Thus, how to develop highly efficient, high-order numerical schemes which preserve the original mass and energy of the LogSE \eqref{RLogSE} is challenging, which motivates this paper.


In this paper, a class of highly efficient, high-order, and mass- and energy-conserving schemes are proposed. The key idea of the scheme include the following two parts: firstly, based on the idea of the supplementary variable method (SVM) proposed by Gong et al.~\cite{GHW2021,GongWH2023}, we reformulate the regularized LogSE model \eqref{RLogSE} into an extended system by introducing two supplementary variables; secondly, we discretize the SVM system by employing a high-order prediction-correction scheme in time and a Fourier pseudo-spectral method in space, respectively. The newly developed scheme has the following advantages:
\begin{itemize}
\item conserve both mass and original energy of the model \eqref{RLogSE};
\item is high-order accuracy in time;
\item only requires solving a constant-coefficient linear systems plus two algebraic equations, which can be solved by the Newton iteration efficiently.

\end{itemize}

The structure of the paper is arranged as follows. In Section \ref{Sec:PM:2}, we reformulate the RlogSE model \eqref{RLogSE} into an equivalent form based on the idea of the supplementary variable method. In Section \ref{Sec:PM:3}, based on the high-order prediction-correction method and the Fourier pseudo-spectral method, a class of high-order mass- and energy-preserving schemes are presented. In Section \ref{Sec:PM:4}, various numerical examples and comparisons are provided to demonstrate the performance of the new scheme. Section \ref{Sec:PM:5} contains a few concluding remarks.

\section{Model reformulation}\label{Sec:PM:2}
In this section, based on the idea of the supplementary variable method\cite{GHW2021,GongWH2023}, we reformulate the RlogSE system \eqref{RLogSE} into a new extended, consistent and well-determined system, which provides an elegant
platform for developing highly efficient, high-order, and mass- and energy-conserving schemes of the RlogSE system \eqref{RLogSE}.

To begin with, we define the inner product in $L^2(\Omega)$ space as
\begin{equation*}
(u,v)=\int_{\Omega}{u}({\bf x})\bar{v}({\bf x}){\rm d}{\bf x}
\end{equation*}
where ${u}({\bf x}),~{v}({\bf x})\in L^2(\Omega)$, $\bar{u}({\bf x})$ is complex conjugate transposition of $u({\bf x})$ and ${\bf x}\in\mathbb{R}^d,\ d=1,2$. Furthermore, the norm in $L^2(\Omega)$ space is denoted as $\|u\|^2=(u,u)$. Then we rewrite the RlogSE system \eqref{RLogSE} into an energy-conserving system
\begin{align}\label{Hamiltonian form}
\partial_tu^{\varepsilon}=\mathcal{G} \frac{\delta \mathcal{E}^{\epsilon}}{\delta \bar{u}^{\varepsilon}},\ \mathcal{G}=- \mathrm{i},
\end{align}
where $\frac{\delta \mathcal{E}^{\epsilon}}{\delta
\bar{u}^{\varepsilon}}$ is the variational derivative of the energy functional $\mathcal{E}^{\epsilon}$ with respect to $\bar{u}^{\varepsilon}$, and $\bar{u}^{\varepsilon}$ is the complex conjugate of $u^{\varepsilon}$. In general, it is challenging task to propose highly efficient, high-order, and mass- and energy-conserving schemes based on the system \eqref{Hamiltonian form}. Motivated by Refs.\cite{GHW2021,GongWH2023}, we first reformulated the original model \eqref{Hamiltonian form} into the following SVM system by introducing two supplementary variables
 \begin{align}\label{SVM-Formulation}
 \left\lbrace
  \begin{aligned}
&\partial_tu^{\epsilon}=\mathcal{G}\frac{\delta
\mathcal E^{\epsilon}}{\delta\bar{u}^{\varepsilon}}+\alpha_{1}g_{1}[u^{\epsilon}]+\alpha _{2}g_{2}[u^{\epsilon}],\\
&\frac {d}{dt}\mathcal{M}^{\epsilon}[u^{\epsilon}]=0,\ \frac {d}{dt}\mathcal{E}^{\epsilon}[u^{\epsilon}]=0,
 \end{aligned}\right.
 \end{align}
 where $\alpha_1$ and $\alpha_2$ are supplementary variables, $g_1[u^{\epsilon}]$ and $g_2[u^{\epsilon}]$ are user-supplied functional of $u^{\epsilon}$. Noticing that the system \eqref{SVM-Formulation} is equivalent to the original system \eqref{Hamiltonian form} as $\alpha_1=0$ and $\alpha_2=0$.

\begin{rmk}
\label{r211}
It is worth noting that there are many flexible ways on how to introduce the supplementary variables method. In following computations, we choose  $g_{1}[u^{\varepsilon}] = \frac{\delta E^{\epsilon}}{\delta
\bar{u}^{\varepsilon}}$ and $g_{2}[{u^{\varepsilon}}] = \frac{\delta M^{\epsilon}}{\delta \bar{u}^{\varepsilon}}$, respectively.
\end{rmk}

\section{High-order mass- and energy-conserving scheme}\label{Sec:PM:3}
In this section, a class of high-order mass- and energy-preserving schemes are proposed for the RlogSE system \eqref{RLogSE} by utilizing the high-order prediction-correction method in time \cite{GongWH2023} and the standard Fourier
pseudo-spectral method in space to discretize the SVM system
\eqref{SVM-Formulation}, respectively.
\subsection{Temporal semi-discretization}\par
In this subsection, we apply a high-order prediction-correction method to the SVM reformulation \eqref{SVM-Formulation} and a new class of temporal semi-discrete schemes are presented, which preserve semi-discrete both mass and energy. Then, we show the proposed scheme can be solved efficiently. Finally, the local existence and uniqueness of solution of supplementary variables and the local truncation errors of the proposed schemes are analyzed based on the idea introduced in \cite{CHMR06,GongWH2023}.

Let $\tau $ be the time step. We define $t_{n} = n \tau$ and $t_{ni} = t_{n} + c_{i} \tau, n = 0, 1,\cdots,$ and denote $u^{\varepsilon,n}$ and $u_{i}^{\varepsilon,n}$ as the numerical approximations of the function $u^{\epsilon}(\textbf{x},t_{n})$ at $t_{n}$ and $t_{ni}$, respectively. Then, a high-order prediction-correction method is employed to discretize the SVM system \eqref{SVM-Formulation} in $[t_{n},t_{n+1} ]$, and one can obtain the following scheme
\begin{shm}\label{svm-scheme1}
	Let $a_{ij}$, $b_{i}$, $c_{i}$ ($i, j=1,2,\cdots,s$) be real numbers and satisfy $c_{i} = \sum_{j=1}^{s} a_{ij}$. For the given $u^{\varepsilon,n}$, we compute $u^{\varepsilon,n+1}$ via the following two steps.
	\begin{enumerate}
		\item \textbf{High-order prediction}: set $u_{i, (0)}^{\varepsilon, n} = u^{\varepsilon, n}$, and compute $u_{i, (m+1)}^{\varepsilon,
n}$,
$\forall$ $i$  from $m=0$ to $K$-1 by solving the following linear system
		\begin{align}\label{semi-prediction-step}
		\left\lbrace
		\begin{aligned}
		&u_{i, (m+1)}^{\varepsilon,n} = u^{\varepsilon,n} + \tau \sum_{j=1}^{s} a_{ij} k_{j, (m+1)}^n, ~i = 2,\cdots,s,\\
		&k_{i, (m+1)}^{n} = \mathrm{i}\Delta u_{i, (m+1)}^{\varepsilon,n}-\mathrm{i}\lambda u_{i, (m)}^{\varepsilon,n}\ln(\varepsilon +|u_{i,
(m)}^{\varepsilon,n}|)^{2}.
		\end{aligned}\right.
		\end{align}
	   Then, we set $u_{i}^{\varepsilon, n\ast}= u_{i, (K)}^{\varepsilon, n}$, $i=1,2,\cdots,s$.
		\item \textbf{High-order correction}: for given $u_{i}^{\varepsilon, n\ast}$, we update $u^{\varepsilon,n+1}$ via
		\begin{align}\label{semi-correction}
        \left\lbrace
        \begin{aligned}
		&u_{i}^{\varepsilon,n} = u^{\varepsilon,n} + \tau \sum_{j=1}^{s} a_{ij} k_{j}^n,\\
        &k_{i}^{n} = \mathrm{i}\Delta u_{i}^{\varepsilon,n}-\mathrm{i}\lambda u_{i}^{\varepsilon,n\ast}\ln(\varepsilon
        +|u_{i}^{\varepsilon,n\ast}|)^{2}+\alpha _{1}^{n}g_{1}[u_{i}^{\varepsilon,n\ast}]+\alpha
        _{2}^{n}g_{2}[u_{i}^{\varepsilon,n\ast}],\\
        &u^{\varepsilon,n+1} = u^{\varepsilon,n} + \tau \sum_{i=1}^{s} b_{i} k_{i}^n,\\
        &\mathcal{M}^{\epsilon}[u^{\varepsilon,n+1}] = \mathcal{M}^{\epsilon}[u^{\varepsilon,n}],\ \mathcal{E}^{\epsilon}[u^{\varepsilon,n+1}] = \mathcal{E}^{\epsilon}[u^{\varepsilon,n}].
        \end{aligned}\right.
		\end{align}	
	\end{enumerate}
\end{shm}

\begin{thm}\label{mass and energy} It is evident that \textbf{Scheme \ref{svm-scheme1}} preserves the following semi-discrete mass and Hamiltonian energy
\begin{align*}
\mathcal{M}^{\epsilon,n}= \mathcal{M}^{\epsilon,0},\ \mathcal{E}^{\epsilon,n}= \mathcal{E}^{\epsilon,0},\ n = 0,1, \cdots ,
\end{align*}
where
\begin{align*}
\mathcal{M}^{\epsilon,n}= (|u^{\epsilon,n}|^{2},1),\
\mathcal{E}^{\epsilon,n}=-(\Delta u^{\epsilon,n},u^{\epsilon,n})+(2\epsilon \lambda |u^{\epsilon,n}|+ 2\lambda(|u^{\epsilon,n}|^{2} - \epsilon^{2}) \ln(\epsilon + |u^{\epsilon,n}|),1).
\end{align*}
\end{thm}
\begin{prf} The conclusion is obvious from the last formula in \textbf{Scheme \ref{svm-scheme1}}. \qed

\end{prf}

Then, we show how to implement \textbf{Scheme \ref{svm-scheme1}} efficiently. To begin with, we denote
\begin{align}
\textbf{k}_{(m+1)}^{n}=(k_{1, (m+1)}^{n},\cdots , k_{s, (m+1)}^{n})^{\top},\ \textbf{u}_{(m)}^{\varepsilon,n}=(u_{1,
(m)}^{\varepsilon,n},\cdots , u_{s, (m)}^{\varepsilon,n})^{\top},
\end{align}
and then rewrite \eqref{semi-prediction-step} as a compact form, as follows:
\begin{align}\label{Operator-matrix}
\textbf{A} \textbf{k}_{(m+1)}^{n}=\textbf{U}(u^{\varepsilon,n}, \textbf{u}_{(m)}^{\varepsilon,n}),
\end{align}
 where $\textbf{A}$ represents the $s\times s$ operator matrix
 \begin{align*}
\textbf{A}= {\left[
\begin{array}{cccc}
1-{\rm i}\tau a_{11}\Delta & -{\rm i}\tau a_{12}\Delta & \cdots & -{\rm i}\tau a_{1s}\Delta \\
-{\rm i}\tau a_{21}\Delta & 1-{\rm i}\tau a_{22}\Delta & \cdots & -{\rm i}\tau a_{2s}\Delta \\
\vdots & \vdots & \ddots &\vdots\\
-{\rm i}\tau a_{s1}\Delta & -{\rm i}\tau a_{s2}\Delta & \cdots & 1-{\rm i}\tau a_{ss}\Delta
\end{array}
\right ]},
\end{align*}
and the components of $\textbf{U}(u^{\varepsilon,n}, \textbf{u}_{(m)}^{\varepsilon,n})$ are given by
\begin{align}
\textbf{U}_{i}(u^{\varepsilon,n}, \textbf{u}_{(m)}^{\varepsilon,n}) = \mathrm{i} \Delta u^{\varepsilon,n}- \mathrm{i} \lambda
u_{i,(m)}^{\varepsilon,n}\ln(\varepsilon +|u_{i,(m)}^{\varepsilon,n}|)^{2} .
\end{align}
Thus, it follows from \eqref{Operator-matrix}
\begin{align}
\textbf{k}_{(m+1)}^{n}=\textbf{A}^{-1} \textbf{U}(u^{\varepsilon,n}, \textbf{u}_{(m)}^{\varepsilon,n}),
\end{align}
where $\textbf{A}^{-1}$ denotes the inverse operator of $\textbf{A}$. Similarly, the system \eqref{semi-correction} is equivalent to
\begin{align}\label {LogRNLS-3.7}
   \textbf{A}\textbf{k}^{n}=\textbf{U}(u^{\varepsilon,n},\textbf{u}^{\varepsilon,n\ast})+\bm{\alpha}_{1}^{n}\textbf{G}_{1}[\textbf{u}^{\varepsilon,n\ast}]+\bm{\alpha}_{2}^{n}\textbf{G}_{2}[\textbf{u}^{\varepsilon,n\ast}],
\end{align}
where ${\bf k}^n=[k_1^n,k_2^n,\cdots,k_s^n]^{\top}$, and the components of $\textbf{G}_{1}[\textbf{u}^{\varepsilon,n\ast}]$ and $\textbf{G}_{2}[\textbf{u}^{\varepsilon,n\ast}]$ are given by
$g_{1}[u_{i}^{\varepsilon,n\ast}]$ and $g_{2}[u_{i}^{\varepsilon,n\ast}]$ respectively. Then, one can deduce from \eqref{LogRNLS-3.7}
\begin{align}
\label {simplified form}
\textbf{k}^{n}=\hat{\textbf{k}}^{n}+\bm{\alpha}_{1}^{n} \textbf{r}_{1}^{n}+\bm{\alpha}_{2}^{n} \textbf{r}_{2}^{n},
\end{align}
where $\hat{\textbf{k}}^{n}= \textbf{A}^{-1}\textbf{U}(u^{\varepsilon,n},\textbf{u}^{\varepsilon,n\ast})$ and $\textbf{r}_{i}^{n}=
\textbf{A}^{-1}\textbf{G}_{i}[\textbf{u}^{\varepsilon,n\ast}], i = 1,2.$ Let $\beta_{1}=\tau \alpha_{1}^{n}$, $\beta_{2}=\tau \alpha_{2}^{n}$, and $\textbf{b}=(b_{1}, b_{2}, \cdots , b_{s})^{\top}$, we obtain from \eqref{simplified form} and the third-last equality of \eqref{semi-correction} that
\begin{align}
\label {LogNLS-002}
u^{\varepsilon,n+1}= \hat{u}^{\varepsilon,n+1}+\beta_{1}
\textbf{R}_{1}^{n}+\beta_{2} \textbf{R}_{2}^{n},\ ,
\end{align}
where $\hat{u}^{\varepsilon,n}= u^{\varepsilon,n} + \tau \textbf{b}^{\top} \hat{\textbf{k}}^{n}$, $\textbf{R}_{1}^{n} = \textbf{b}^{\top} \textbf{r}_{1}^{n}$ and $\textbf{R}_{2}^{n} = \textbf{b}^{\top} \textbf{r}_{2}^{n}$. Finally, plugging \eqref{LogNLS-002} into the last equality of \eqref{semi-correction} yields two scalar nonlinear algebraic equations for ${\bf \beta}=(\beta_1,\beta_2)^{\top}$, as follows:
\begin{align}
\label {algebraic equation1}
&\mathcal E^{\epsilon}[\hat{u}^{\varepsilon,n}+\beta_{1} \textbf{R}_{1}^{n}+\beta_{2} \textbf{R}_{2}^{n}] -
\mathcal E^{\epsilon}[u^{\epsilon,n}]=0,\hfill\hfill\\
\label {algebraic equation2}
&\mathcal M^{\epsilon}[\hat{u}^{\varepsilon,n}+\beta_{1} \textbf{R}_{1}^{n}+\beta_{2} \textbf{R}_{2}^{n}]-
\mathcal M^{\epsilon}[u^{\epsilon,n}]=0,
\end{align}
which can be solved efficiently using the Newton iteration with the initial guess of $(0,0)^{\top}$.

\begin{thm}\label{RK order}\cite{GongWH2023}~If the underlying RK method is of order $p$, then the prediction-correction RK scheme is of
order $min \{ p,K \}$, i.e.,
\begin{align}
u_{(K)}^{\varepsilon,n+1} - u^{\varepsilon}(t_{n+1}) = \mathcal{O}(\tau^{min\{p,K\}+1}).
\end{align}
\end{thm}

Based on the ideas presented in \cite{CHMR06,GJZ2024,GongWH2023}, we show the existence and uniqueness of ${\bf \beta}$ in \eqref{algebraic equation1}-\eqref{algebraic equation2},  and investigate the order in local truncation error.
\begin{thm} Let
\begin{align*}
\mathcal{J}=
\left[ {\begin{array}{cc}
\Re(\frac {\delta \mathcal E^{\epsilon}}{\delta\bar{u}^{\varepsilon}}[u^{\epsilon,n}],g_{1}[u^{\epsilon, n}]) & \Re(\frac {\delta
\mathcal E^{\epsilon}}{\delta \bar{u}^{\varepsilon}}[u^{\epsilon,n}],g_{2}[u^{\epsilon, n}])\\
\Re(\frac {\delta M^{\epsilon}}{\delta \bar{u}^{\varepsilon}}[u^{\epsilon,n}],g_{1}[u^{\epsilon, n}]) & \Re(\frac {\delta
\mathcal M^{\epsilon}}{\delta \bar{u}^{\varepsilon}}[u^{\epsilon,n}],g_{2}[u^{\epsilon, n}])\\
\end{array} }\right ],
\end{align*}
and suppose $|\mathcal{J} | \neq 0$, there exists a $\tau^{\ast}> 0$ such that the equations \eqref{algebraic equation1} and \eqref{algebraic equation2} defines a unique function ${\bf \beta}={\bf\beta}(\tau)$ for all $\tau \in [0, \tau^{\ast}]$.
\end{thm}
\begin{prf}
To begin with, we define the function $\textbf{\textit{F}}(\tau,
{\bf\beta})=(F_{1}(\tau,\beta),F_{2}(\tau, \beta))^{\top}$, as follows: 
\begin{align}
&F_{1}(\tau, \beta) = \mathcal E^{\epsilon}[\hat{u}^{\varepsilon,n}+\beta_{1} \textbf{R}_{1}^{n}+\beta_{2} \textbf{R}_{2}^{n}] -
\mathcal E^{\epsilon}[u^{\epsilon,n}],\hfill\hfill\\
&F_{2}(\tau, \beta) = \mathcal M^{\epsilon}[\hat{u}^{\varepsilon,n}+\beta_{1} \textbf{R}_{1}^{n}+\beta_{2} \textbf{R}_{2}^{n}]  -
\mathcal M^{\epsilon}[u^{\epsilon,n}].
\end{align}
Considering the smoothness of $\textbf{\textit{F}}$ and
\begin{align}
\textbf{\textit{F}}(0, \textbf{0}) = \textbf{0}, \ \left|\frac{\partial}{\partial {\beta} }\textbf{\textit{F}}(0,
\textbf{0})\right| = |\mathcal J| \neq 0,
\end{align}
according to the implicit function theorem, there exists a $\tau^{\ast} > 0$ such that equations $\textbf{\textit{F}}(\tau,
{\beta})  = \textbf{0}$ defines the unique smooth functions $\beta_{1} = \beta_{1} (\tau)$ and $\beta_{2} = \beta_{2} (\tau)$ satisfying
${\bf \beta}(\textbf{0}) = \textbf{0}$ and $\textbf{\textit{F}}(\tau, {\bf \beta}(\tau)) = \textbf{0}$ for all $\tau \in [0,
\tau^{\ast}]$.\qed
\end{prf}
\begin{thm}
If the underlying RK method is of order $p$, then the proposed \textbf{Scheme \ref{svm-scheme1}} is of order $\hat{p}$, where $\hat{p} = min \{ p,K+1\}$.
\end{thm}
\begin{prf}
Notice that $\hat{u}^{\varepsilon,n+1}$ satisfies
\begin{align}\label{proof}
\left\lbrace
\begin{aligned}
&u_{i,(K+1)}^{\varepsilon,n} = u^{\varepsilon,n} + \tau \sum_{j=1}^{s} a_{ij} k_{j,(K+1)}^n,\\
&k_{i,(K+1)}^{n} = \mathrm{i}\Delta u_{i,(K+1)}^{\varepsilon,n}-\mathrm{i}\lambda u_{i,(K)}^{\varepsilon,n}\ln(\varepsilon
+|u_{i,(K)}^{\varepsilon,n}|)^{2},\\
&\hat{u}^{\varepsilon,n+1} = u^{\varepsilon,n} + \tau \sum_{i=1}^{s} b_{i} k_{i,(K+1)}^n.
\end{aligned}\right.
\end{align}
According to $\textbf{Theorem \ref{RK order}}$, we obtain
\begin{align}
\hat{u}^{\varepsilon,n+1} - u^{\varepsilon}(t_{n+1}) = \mathcal{O}(\tau^{min\{p,K+1\}+1}).
\end{align}
By the Taylor expansion, we then obtain
\begin{align}\label{LogNLS-3.18}
&\mathcal E^{\epsilon}[\hat{u}^{\varepsilon,n+1}] = \mathcal E^{\epsilon}[u^{\varepsilon}(t_{n+1})]+\mathcal{O}(\tau^{min\{p,K+1\}+1}),\\\label{LogNLS-3.19}
&\mathcal M^{\epsilon}[\hat{u}^{\varepsilon,n+1}]= \mathcal M^{\epsilon}[u^{\varepsilon}(t_{n+1})]+\mathcal{O}(\tau^{min\{p,K+1\}+1}),
\end{align}
Denote
\begin{align}
 \mathcal E_{(K)}^{\epsilon,n+1} = \mathcal E^{\epsilon}[u^{\varepsilon,n}], \ \mathcal M_{(K)}^{\epsilon,n+1} = \mathcal M^{\epsilon}[u^{\varepsilon,n}].
\end{align}
we obtain from \eqref{LogNLS-3.18} and \eqref{LogNLS-3.19}
\begin{align}
&\mathcal E^{\epsilon}[\hat{u}^{\varepsilon,n+1}] - \mathcal E_{(K)}^{\epsilon, n+1}= \mathcal E^{\epsilon}[\hat{u}^{\varepsilon,n+1}] - \mathcal E^{\epsilon}[u^{\varepsilon}(t_n)]= \mathcal{O}(\tau^{min\{p,K+1\}+1}),\\
&\mathcal M^{\epsilon}[\hat{u}^{\varepsilon,n+1}] - \mathcal M_{(K)}^{\epsilon, n+1} =\mathcal M^{\epsilon}[\hat{u}^{\varepsilon,n+1}] - \mathcal M[u^{\varepsilon}(t_n)] = \mathcal{O}(\tau^{min\{p,K+1\}+1}).
\end{align}
Then, we expand $\textbf{\textit{F}}(\tau, {\beta})$ at ${\bf \beta}= (0, 0)^{\top}$ to obtain
\begin{align}
\textbf{\textit{F}}(\tau, {\bf \beta})=\textbf{\textit{F}}(\tau, \textbf{0})+ \frac{\partial}{\partial {\bf \beta}
}\textbf{\textit{F}}(\tau, \textbf{0}){\bf \beta} + \mathcal{O}(\beta_{1}^{2}+\beta_{2}^{2}),
\end{align}
where
\begin{align*}
\textbf{\textit{F}}(\tau, \textbf{0})= \Big[ \mathcal  E^{\epsilon}[\hat{u}^{\varepsilon,n+1}] - \mathcal  E_{(K)}^{\epsilon, n+1}, \mathcal  M^{\epsilon}[\hat{u}^{\varepsilon,n+1}] - \mathcal  M_{(K)}^{\epsilon, n+1}\Big]^{\top} = \mathcal{O}(\tau^{min\{p,K+1\}+1}),
\end{align*}
and
\begin{align*}
\frac{\partial}{\partial {\bf \beta}}\textbf{\textit{F}}(\tau, \textbf{0}) = \frac{\partial}{\partial {\bf \beta}}
\textbf{\textit{F}}(0, \textbf{0})+ \mathcal{O}(\tau) = \mathcal {J} + \mathcal{O}(\tau).
\end{align*}
This implies
\begin{align}
\beta_{1} = \beta_{1} (\tau) = \mathcal{O}(\tau^{\hat{p}+1}),\ \beta_{2} = \beta_{2}(\tau) = \mathcal{O}(\tau^{\hat{p}+1}),\
\hat{p}=min \{ p,K+1 \}.
\end{align}
Thus, the proposed scheme is of order $\hat{p}$.\qed
\end{prf}


\begin{rmk}\label {rk-Gauss}
In the following calculations, we mainly focus on the 2-stage Gauss RK method, where the coefficients is displayed as follows \cite{HLWbook2006}:
\begin{table}[H]
\caption{~Coefficients of Gauss RK methods of order 4.}\label{Tab:Gauss-cllocation-method}
\centering
\begin{tabular}{c|cc}
$\frac{1}{2}-\frac{\sqrt{3}}{6}$ &$\frac{1}{4}$ & $\frac{1}{4}- \frac{\sqrt{3}}{6}$\\
$\frac{1}{2}+\frac{\sqrt{3}}{6}$ &$\frac{1}{4}+ \frac{\sqrt{3}}{6}$ &$\frac{1}{4}$ \\
\hline
                                 &$\frac{1}{2}$&$\frac{1}{2}$
\end{tabular}
\end{table}
\end{rmk}

\subsection{Full discretization}\par
In this subsection, we employ the standard Fourier pseudo-spectral method to discretize the semi-discrete \textbf{Scheme \ref{svm-scheme1}} as the periodic boundary condition is considered.

To begin with, we let $\Omega=[x_{L},x_{R}]\times[y_{L},y_{R}]$, and choose the mesh sizes $h_x=(x_{R}-x_{L})/N_x$ and $h_y=(y_{R}-y_{L})/N_y$ with two even positive integers $N_x$ and $N_y$; denote $u_{j,k}^{\epsilon,n}$ and $(u_{j,k})_{i}^{\epsilon,n}$ as the numerical approximations of  $u^{\epsilon}(x_j,y_k,t_n)$ and $u^{\epsilon}(x_j,y_k,t_{ni})$ for $j=0,1,\cdots,N_x,\ k=0,1,2,\cdots,N_y,\ n=0,1,\cdots,$ respectively. Then, we let
  \begin{equation*}
  {\bf U}^{\epsilon,n}=(u_{j,k}^{\epsilon,n})_{N_x\times N_y}=\left(\begin{array}{ccccc}
 u_{0,0}^{\epsilon,n}&u_{0,1}^{\epsilon,n} &\cdots&u_{0,N_y-1}^{\epsilon,n}\vspace{1mm}\\
u_{1,0}^{\epsilon,n}&u_{1,1}^{\epsilon,n} &\cdots&u_{1,N_y-1}^{\epsilon,n}\vspace{1mm}\\
\vdots&\vdots &\ddots&\vdots\vspace{1mm}\\
u_{N_x-1,0}^{\epsilon,n}&,u_{N_x-1,1}^{\epsilon,n} &,\cdots&,u_{N_x-1,N_y-1}^{\epsilon,n}\vspace{1mm}\\
\end{array} \right)_{Nx\times N_y}.
  \end{equation*}
  be the grid matrix function. Subsequently, for any two grid matrix functions ${\bf U}^n$ and ${\bf V}^n$, we define the discrete $L^2$-inner product and norm as, respectively,
\begin{equation*}
\langle {\bf U}^{\epsilon,n},{\bf V}^{\epsilon,n}\rangle _{h}=h_x h_y\sum_{j=0}^{N_x-1}\sum_{k=0}^{N_y-1} u_{j,k}^{\epsilon,n}\bar{v}_{j,k}^{\epsilon,n},\
\|{\bf U}^{\epsilon,n}\|_{h}^2=\langle {\bf U}^{\epsilon,n},{\bf U}^{\epsilon,n}\rangle_{h}.
\end{equation*}
Finally, we define another operator $``\odot"$ for element by element multiplication between two matrix functions of same sizes as
\begin{equation*}
({\bf U}^{\epsilon,n}\odot{\bf V}^{\epsilon,n})_{j,k}=({\bf V}^{\epsilon,n}\odot{\bf U}^{\epsilon,n})_{j,k}=(u_{j,k}^{\epsilon,n}\cdot v_{j,k}^{\epsilon,n}).
\end{equation*}
Then, consider the interpolation space
\begin{align*}
&\mathcal T_N=\text{span}\{X_{j}(x)Y_{k}(y),\ 0\leq j\leq N_x-1, 0\leq k\leq N_y-1\}
\end{align*}
where $X_{j}(x)$ and $Y_{k}(y)$ are trigonometric polynomials of degree $N_{x}/2$ and $N_{y}/2$,
given, respectively, by
\begin{align*}
  &X_{j}(x)=\frac{1}{N_{x}}\sum_{l=-N_{x}/2}^{N_{x}/2}\frac{1}{a_{l}}e^{\text{i}l\mu_{x} (x-x_{j})},\ a_{l}=\left \{
 \aligned
 &1,\ |l|<\frac{N_x}{2},\\
 &2,\ |l|=\frac{N_x}{2},
 \endaligned
 \right.\ \mu_{x}=\frac{2\pi}{b-a},\\
 & Y_{k}(y)=\frac{1}{N_{y}}\sum_{l=-N_{y}/2}^{N_{y}/2}\frac{1}{b_{l}}e^{\text{i}l\mu_{y} (y-y_{k})},\ b_{l}=\left \{
 \aligned
 &1,\ |l|<\frac{N_y}{2},\\
 &2,\ |l|=\frac{N_y}{2}.
 \endaligned
 \right.\ \mu_{y}=\frac{2\pi}{d-c}.
\end{align*}
We define interpolation operator $\mathcal I_N: C(\Omega)\rightarrow \mathcal T_N$:
\begin{align}\label{Rlog-NLS-IO}
\mathcal I_Nu^{\epsilon}(x,y,t)=\sum\limits_{j=0}^{N_x-1}\sum\limits_{k=0}^{N_y-1}
u_{j,k}^{\epsilon} X_j(x)Y_k(y),
\end{align}
where $u_{j,k}^{\epsilon}:=u^{\epsilon}(x_j,y_k,t)$. We compute partial derivatives in $x$ and $y$, respectively at the collocation points ($x_{j},y_{k}$) yields
\begin{align*}
\partial^{s_x}_x\partial^{s_y}_y \mathcal I_{N}u^{\epsilon}(x_j, y_k)=\sum\limits_{p=0}^{N_x-1}\sum\limits_{q=0}^{N_y-1}
u_{p,q}^{\epsilon}(\textbf{D}^{x}_{s_x})_{j, p}(\textbf{D}^{y}_{s_y})_{k, q}=\big(\textbf{D}^{x}_{s_x}\textcircled{x}\textbf{D}^{y}_{s_y}\textcircled{y}{\bf U }^{\epsilon}\big)_{j,k},
\end{align*}
where $\textbf{D}^{x}_{s_x}$ and $\textbf{D}^{y}_{s_y}$ are $N_x \times N_x$ and $N_y \times N_y$ matrices, respectively, with
entries given by
\begin{align*}
(\textbf{D}^{x}_{s_x})_{j, p}=\frac{d^{s_x}X_p(x_j)}{dx^{s_x}}, ~~(\textbf{D}^{y}_{s_y})_{k, q}=\frac{d^{s_y}Y_q(y_k)}{dy^{s_y}}.
\end{align*}
\begin{rmk}\label {r311}
It is worth noting that the $\textbf{D}^{a}_{s_a}$ $(a=x ~or~ y)$ can be expressed as follows \cite{GCWCiCP2014,STW2011}:
\begin{align*}
\textbf{D}^{a}_{s_a}=\left\{
\begin{aligned}
\mathcal{F}_{N_a}^{-1} {\Lambda}^{s_a}\mathcal{F}_{N_a},~~s_a\mbox{ is an odd number},\\
\mathcal{F}_{N_a}^{-1} \widetilde{\Lambda}^{s_a} \mathcal{F}_{N_a},~~s_a\mbox{ is an even number},\\
\end{aligned}
\right.
\end{align*}
where
\begin{align*}
&{\Lambda}^{s_a}={\rm diag}\big ({\rm i}\mu_a[0,1,\cdots,\dfrac{N_a}{2}-1,0,-\dfrac{N_a}{2}+1,\cdots,-2,-1]\big ),\\
&\widetilde{\Lambda}^{s_a}={\rm diag}\big ({\rm i}\mu_a[0,1,\cdots,\dfrac{N_a}{2}-1,\dfrac{N_a}{2},-\dfrac{N_a}{2}+1,\cdots,-2,-1]\big ),
\end{align*}
and $\mathcal{F}_{N_a}$ is the discrete Fourier transform (DFT) and $\mathcal{F}_{N_a}^{-1}$ represents the inverse discrete Fourier transform.
\end{rmk}
Then, applying the Fourier pseudo-spectral method described as above to \textbf{Scheme \ref{svm-scheme1}}, we obtain a fully discrete scheme of the system \eqref{SVM-Formulation}, as follows:

\begin{shm}\label{scheme2}
	Let $a_{ij}$, $b_{i}$, $c_{i}$ ($i, j=1,2,...,s$) be real numbers and satisfy $c_{i} = \sum_{j=1}^{s} a_{ij}$.  For given $u^{\varepsilon,n}$, an s-stage RK Fourier pseudo-spectral method is given by
	\begin{enumerate}
		\item \textbf{High-order prediction}: set ${\bf U}_{i, (0)}^{\varepsilon, n} = {\bf U}^{\varepsilon, n}$, and compute ${\bf U}_{i, (m+1)}^{\varepsilon,n}$,
$\forall$ $i$  from $m=0$ to $K$-1 by solving the following linear system
		\begin{align}\label{prediction-step}
		\left\lbrace
		\begin{aligned}
		&{\bf U}_{i, (m+1)}^{\varepsilon,n} = {\bf U}^{\varepsilon,n} + \tau \sum_{j=1}^{s} a_{ij} k_{j, (m+1)}^n, ~i = 2,\cdots,s,\\
		&k_{i, (m+1)}^{n} = \mathrm{i}\Delta_{h} {\bf U}_{i, (m+1)}^{\varepsilon,n}-\mathrm{i}\lambda {\bf U}_{i, (m)}^{\varepsilon,n}\odot\ln(\varepsilon +|{\bf U}_{i,
(m)}^{\varepsilon,n}|)^{2},
		\end{aligned}\right.
		\end{align}
	   where $\Delta_{h} {\bf U}^{\varepsilon,n} = \bf{D}^{x}_{2} \textcircled{x}{\bf U}^{\varepsilon,n} + \bf{D}^{y}_{2}\textcircled{y}{\bf U}^{\varepsilon,n}$. Then, we set ${\bf U}_{i}^{\varepsilon, n\ast}= {\bf U}_{i, (M)}^{\varepsilon, n}$, $i=1,2,\cdots,s$.
		\item \textbf{ High-order correction}: for given ${\bf U}_{i}^{\varepsilon, n\ast}$, we update ${\bf U}^{\varepsilon,n+1}$ via
		\begin{align}\label{correction}
        \left\lbrace
        \begin{aligned}
		&{\bf U}_{i}^{\varepsilon,n} = {\bf U}^{\varepsilon,n} + \tau \sum_{j=1}^{s} a_{ij} k_{j}^n,\\
        &k_{i}^{n} = \mathrm{i}\Delta_{h} {\bf U}_{i}^{\varepsilon,n}-\mathrm{i}\lambda {\bf U}_{i}^{\varepsilon,n\ast}\odot\ln(\varepsilon
        +|{\bf U}_{i}^{\varepsilon,n\ast}|)^{2}+\alpha _{1}^{n}g_{1}[{\bf U}_{i}^{\varepsilon,n\ast}]+\alpha
        _{2}^{n}g_{2}[{\bf U}_{i}^{\varepsilon,n\ast}],\\
        &{\bf U}^{\varepsilon,n+1} = {\bf U}^{\varepsilon,n} + \tau \sum_{i=1}^{s} b_{i} k_{i}^n,\\
        &\mathcal M_{h}^{\epsilon}[{\bf U}^{\varepsilon,n+1}] = \mathcal M_{h}^{\epsilon}[{\bf U}^{\varepsilon,n}],\ \mathcal E_{h}^{\epsilon}[{\bf U}^{\varepsilon,n+1}] = \mathcal E_{h}^{\epsilon}[{\bf U}^{\varepsilon,n}].
        \end{aligned}\right.
		\end{align}	
	\end{enumerate}
\end{shm}

\begin{thm}\label{full discrete mass and energy} \textbf{Scheme \ref{scheme2}} satisfies the following discrete mass and energy conservation laws
\begin{align*}
\mathcal M_{h}^{\epsilon,n}= \mathcal M_{h}^{\epsilon,0},\ \mathcal E_{h}^{\epsilon,n}= \mathcal E_{h}^{\epsilon,0},\ n = 0,1, \cdots,
\end{align*}
where
\begin{align*}
\mathcal M^{\epsilon,n}= \|{\bf U}^{\epsilon,n}\|_{h}^2,\ \mathcal E_{h}^{\epsilon,n}=-\langle \Delta_{h} {\bf U}^{\epsilon,n},{\bf U}^{\epsilon,n}\rangle_{h}+\langle 2\epsilon \lambda |{\bf U}^{\epsilon,n}|+ 2\lambda (|{\bf U}^{\epsilon,n}|^{2}- \epsilon^{2})\odot\ln(\epsilon + |{\bf U}^{\epsilon,n}|),{\bf 1}\rangle_{h}.
\end{align*}
\end{thm}
\begin{prf} The conclusion is obvious from the last formula in \textbf{Scheme \ref{scheme2}}. \qed

\end{prf}

\section{Numerical experiments}\label{Sec:PM:4}
In this section, the convergence, accuracy and conservative properties of the proposed scheme are verified through the presentation of several numerical results. For brevity, in the rest of this paper, the fourth-order SVM method (abbreviated as {\bf SVM4}) are only used for demonstration purposes. Additionally, we compare it with the fourth-order IEQ method described in Ref.~\cite{qian2023} (abbreviated as {\bf IEQ4}) where the diagonally implicit Runge-Kutta method is replaced by the Gauss method of order 4 (see Table \ref{Tab:Gauss-cllocation-method}).

To quantify the numerical errors, we introduce the $L^2$-error function and convergence order as, respectively
\begin{align*}
 e^{\epsilon,\tau,h}_{2}(t_{n}) = \| {\bf U}^{\epsilon}(\cdot ,t_{n}) - {\bf U}^{\epsilon,n} \|_{h},\ {\rm Order}=\log_2\Bigg(e_{2}^{\epsilon,2\tau,h}(t_n)/e_{2}^{\epsilon,\tau,h}(t_n)\Bigg).
\end{align*}
Furthermore, we also define the relative residual functions on the mass and energy as, respectively
\begin{align*}
e_{\mathcal M}(t_n)=\left|\mathcal{M}_{h}^{\epsilon,n}-\mathcal{M}_{h}^{\epsilon,0}\right|/|\mathcal{M}_{h}^{\epsilon,0}|,  \ e_{\mathcal E}(t_n)=\left|\mathcal{E}_{h}^{\epsilon,n}-\mathcal{E}_{h}^{\epsilon,0}\right|/|\mathcal{E}_{h}^{\epsilon,0}|.
\end{align*}


\subsection{RlogSE in 1D}\label{sec:4.1}

\begin{ex} [\textbf{Accuracy confirmation in 1D}]\label{exmp1-accuracy} In this example, we will test temporal numerical error, convergence order of the \textbf{SVM} and \textbf{IEQ4} scheme for the wave function ${u}^{\epsilon}(\cdot ,T=1)$ of the RlogSE \eqref{RLogSE} in 1D by taking $x_{L}=-x_{R}=-16$, the Fourier node 512, the parameters $\lambda = -1,\ \epsilon = 1.0\ \times 10^{-15}$ and the following initial condition
\begin{align}\label{EX1-1D}
u_{0}(x) = \sqrt{-\lambda / \pi} e^{ix + \frac{\lambda}{2} x^{2}}, \ x \in \Omega=[x_{L},x_{R}].
\end{align}
Due to the exact solution is not known, we take the numerical solution produced by the proposed {\bf SVM4} with the time step $\tau = 1\times 10^{-4}$  as a ``reference solution".
\end{ex}

\begin{table}[H]
\tabcolsep=9pt
\footnotesize
\renewcommand\arraystretch{1.1}
\centering
\caption{Temporal numerical error and convergence order at $T = 1$ in example \ref{exmp1-accuracy}}
\label{Tab:EX1-1D-order}
\begin{tabular*}{\textwidth}[h]{@{\extracolsep{\fill}}c c c c c c c}\hline
{Scheme\ \ } &{}&{$\tau_0=1/40$} &{$\tau_0/2$} &{$\tau_0/4$} &{$\tau_0/8$} &{$\tau_0/16$} \\
\hline
\multirow{2}{*}{SVM4}
  &{$e^{\epsilon,\tau,h}_{2}$}&{6.70e-06}&{4.19e-07}&{2.62e-08}&{1.64e-09}&{1.02e-10}
  \\[1ex]
 {}  &{Order}  &{-}& {3.998}&{3.999}&{4.000}&{4.000} \\[1ex]
\multirow{2}{*}{IEQ4}
   &{$e^{\epsilon,\tau,h}_{2}$}&{6.79e-06}&{4.25e-07}&{2.66e-08}&{1.66e-09}&{1.04e-10}
   \\[1ex]
{}  &{Order}  &{-}& {3.998}&{3.999}&{4.000}&{4.000}\\\hline
\end{tabular*}
\end{table}


 Table \ref{Tab:EX1-1D-order} shows the numerical errors and convergence orders of the {\bf SVM4} and {\bf IEQ4} scheme. From the Table, it is clear to observe that the two schemes are fourth order accurate in time, and the errors produced by {\bf SVM4} are much smaller than the ones produced by {\bf IEQ4}.


\begin{ex}\label{EX4-2}[\textbf{Long-time evolution of dynamical behaviors and conservation laws in 1D}]
In this example, we will employ the \textbf{SVM4} scheme to investigate the long time dynamical behaviors and conservation laws of the RlogSE \eqref{RLogSE} in 1D  by choosing the Fourier node 1024 and the time step $\tau = 5 \times 10^{-3}$ with the parameters $\lambda = -1,\ \epsilon = 1.0\ \times 10^{-15}$ and the following initial condition~\cite{qian2023}
\begin{align}
\label {EX2-1D}
u_{0}(x)=\sum _{k=1}^{2}b_{k}e^{-\frac{a_{k}}{2}(x-x_{k})^{2}+iv_{k}x},\ x\in \Omega=[x_{L}, x_{R}],
\end{align}
where $a_{k}, b_{k}, x_{k}$ and $v_{k},\ k=1,2$  are real constants which are chosen as follows:
\begin{itemize}
\item Case I: $x_{L}=-x_{R}=16, x_{1}=-x_{2}=-5, v_{k}=0, a_{k}=b_{k}=1, k=1, 2$;
\item Case II: $x_{L}=-x_{R}=-40, x_{1}=-x_{2}=-3, v_{k}=0, a_{k}=b_{k}=1, k=1, 2$;
\item Case III: $x_{L}=-x_{R}=-50, x_{1}=-x_{2}=-30, v_{1}=-v_{2}=2, a_{k}=b_{k}=1, k=1, 2$;
\item Case IV: $x_{L}=-x_{R}=-50, x_{1}=-x_{2}=-30, v_{1}=-v_{2}=15, a_{k}=b_{k}=1, k=1, 2$.
\end{itemize}
\end{ex}

\begin{figure}[H]
\centering\begin{minipage}[t]{65mm}
\includegraphics[width=65mm]{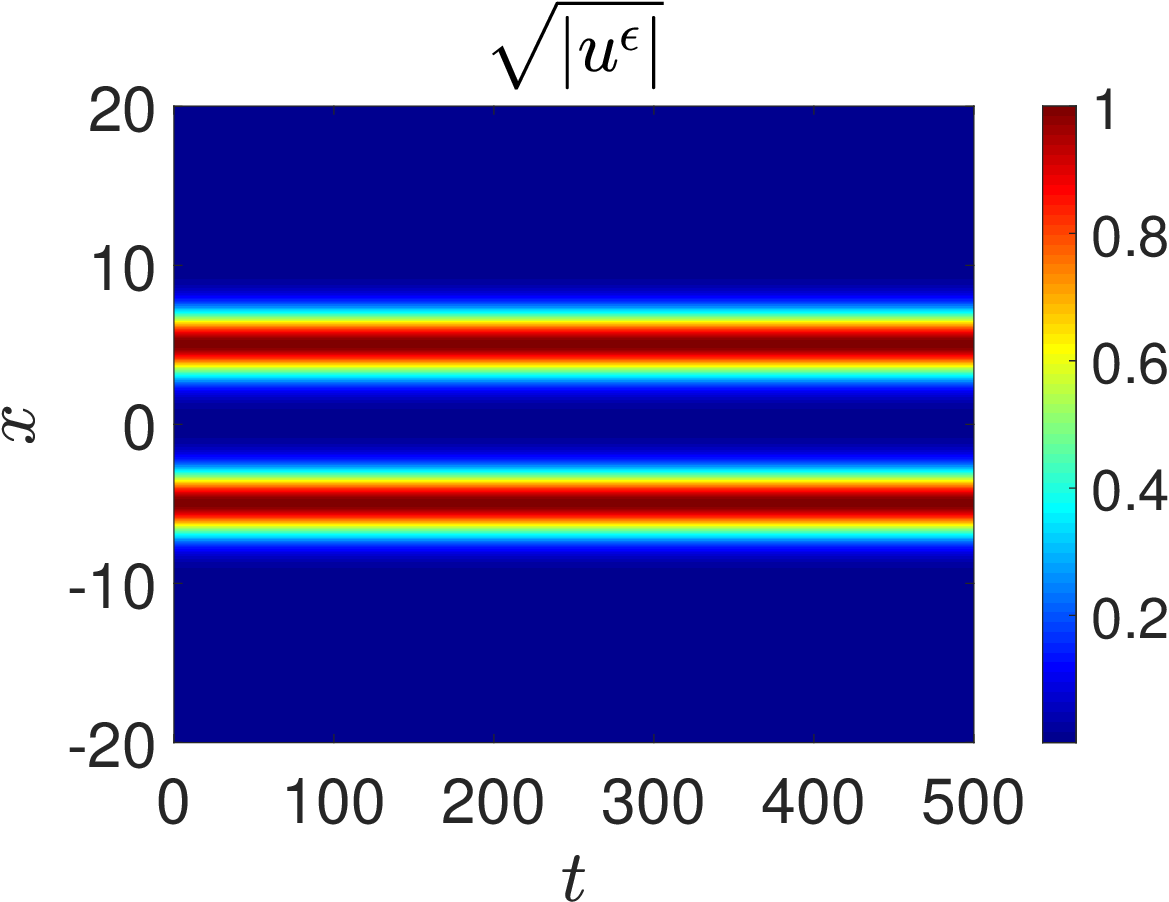}
\end{minipage}
\hspace{0.03\textwidth} 
\begin{minipage}[t]{65mm}
\includegraphics[width=65mm]{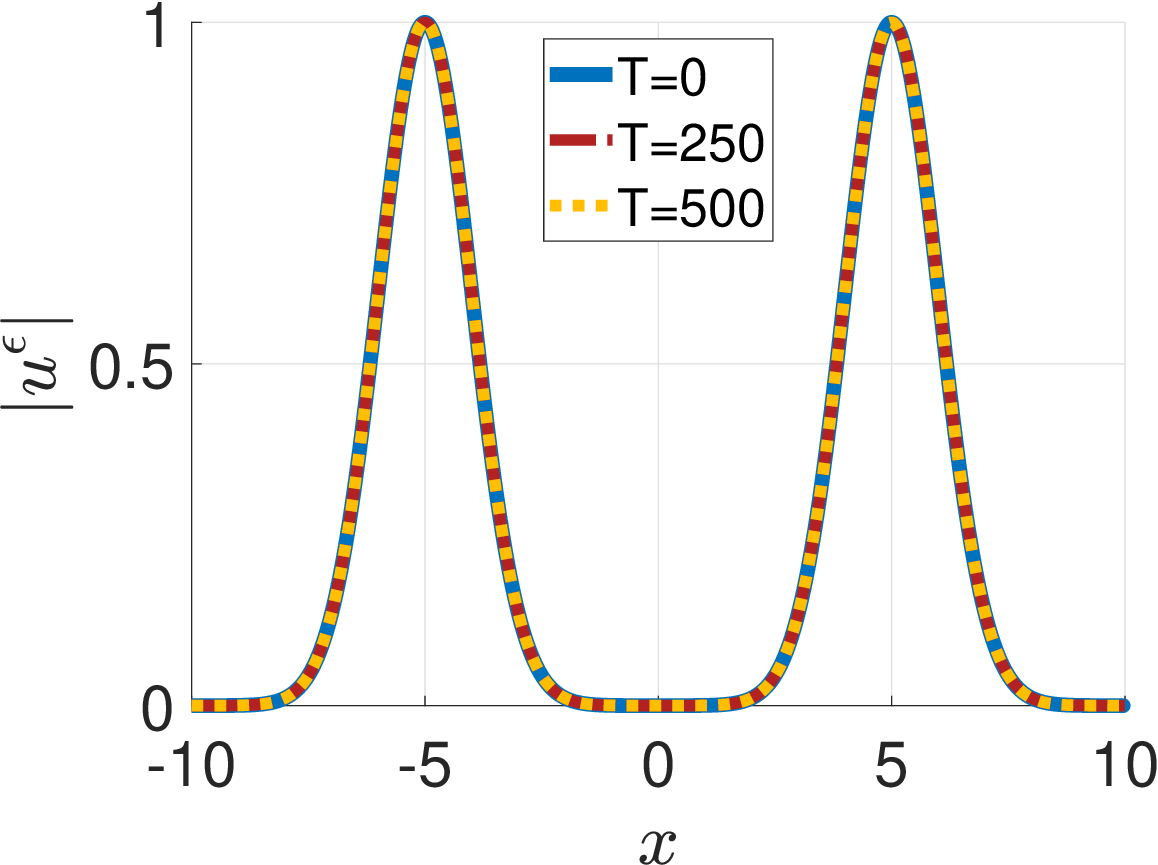}
\end{minipage}
\caption{ Case I : Evolution of $\sqrt{\left| u^\epsilon \right|}$ from $t = 0$ to $t = 500$ (left), and plots of $\left| u^\epsilon \right|$ at different times (right) produced by {\bf SVM4} in example \ref{EX4-2}.}\label{1d-scheme:hotcase1}
\end{figure}

\begin{figure}[H]
\centering\begin{minipage}[t]{65mm}
\includegraphics[width=65mm]{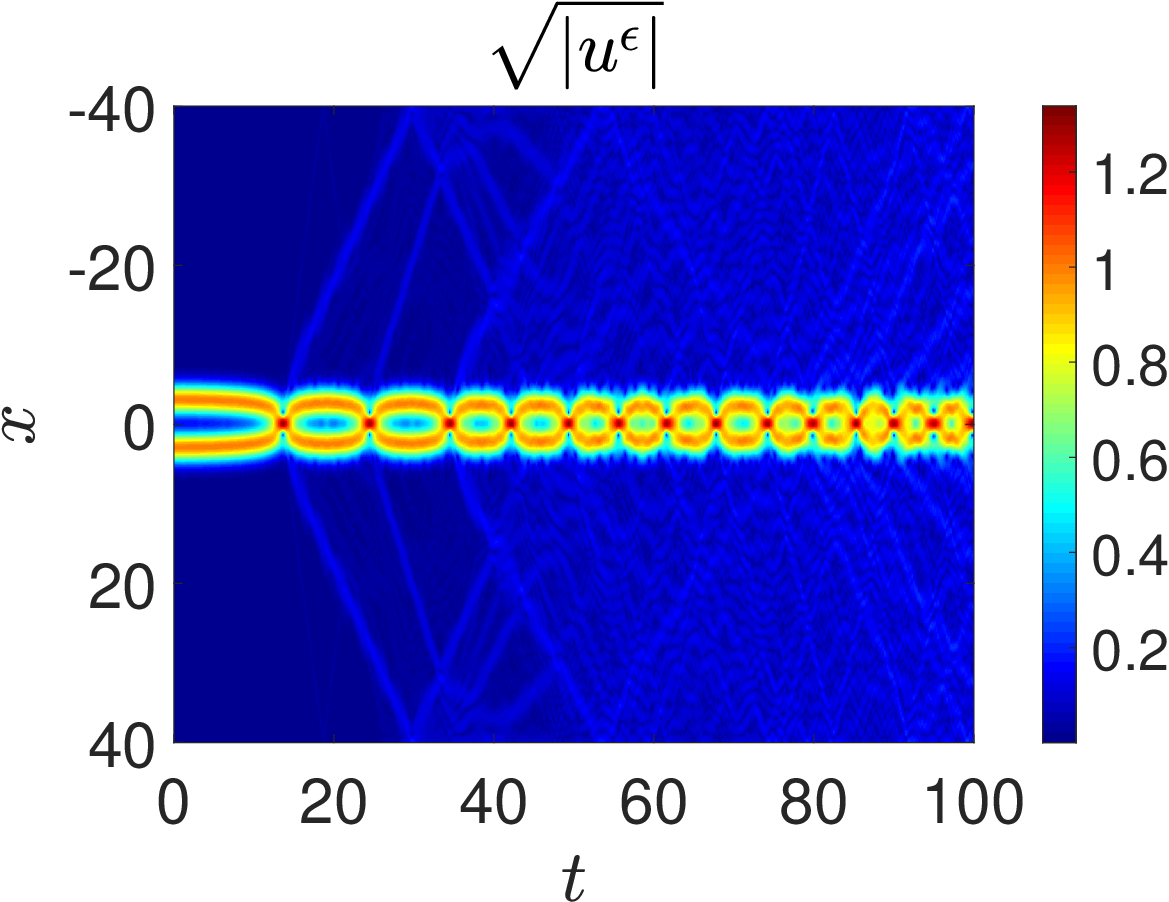}
\end{minipage}
\hspace{0.03\textwidth}
\begin{minipage}[t]{65mm}
\includegraphics[width=65mm]{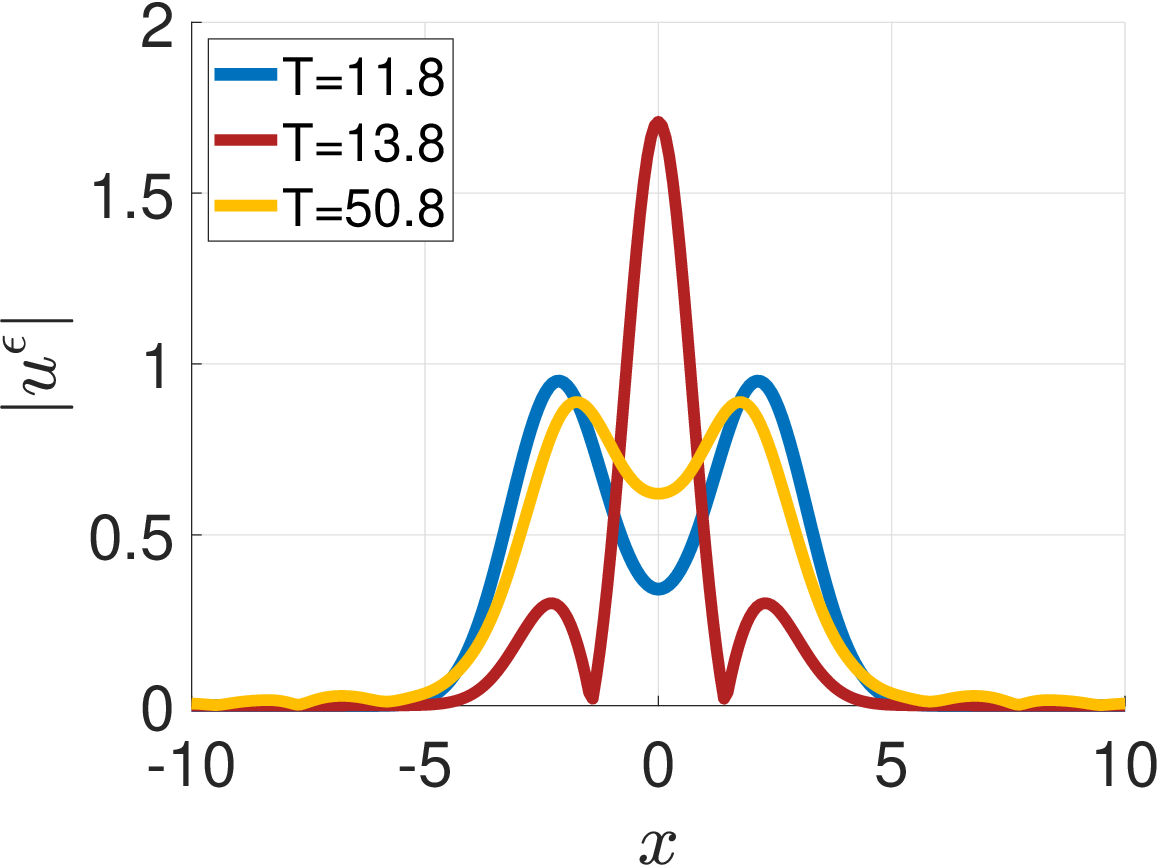}
\end{minipage}
\caption{Case II : Evolution of $\sqrt{\left| u^\epsilon \right|}$ from $t = 0$ to $t = 100$ (left), and plots of $\left| u^\epsilon \right|$ at different times (right) produced by {\bf SVM4} in example \ref{EX4-2}.}\label{1d-scheme:hotcase2}
\end{figure}

\begin{figure}[H]
\centering\begin{minipage}[t]{65mm}
\includegraphics[width=65mm]{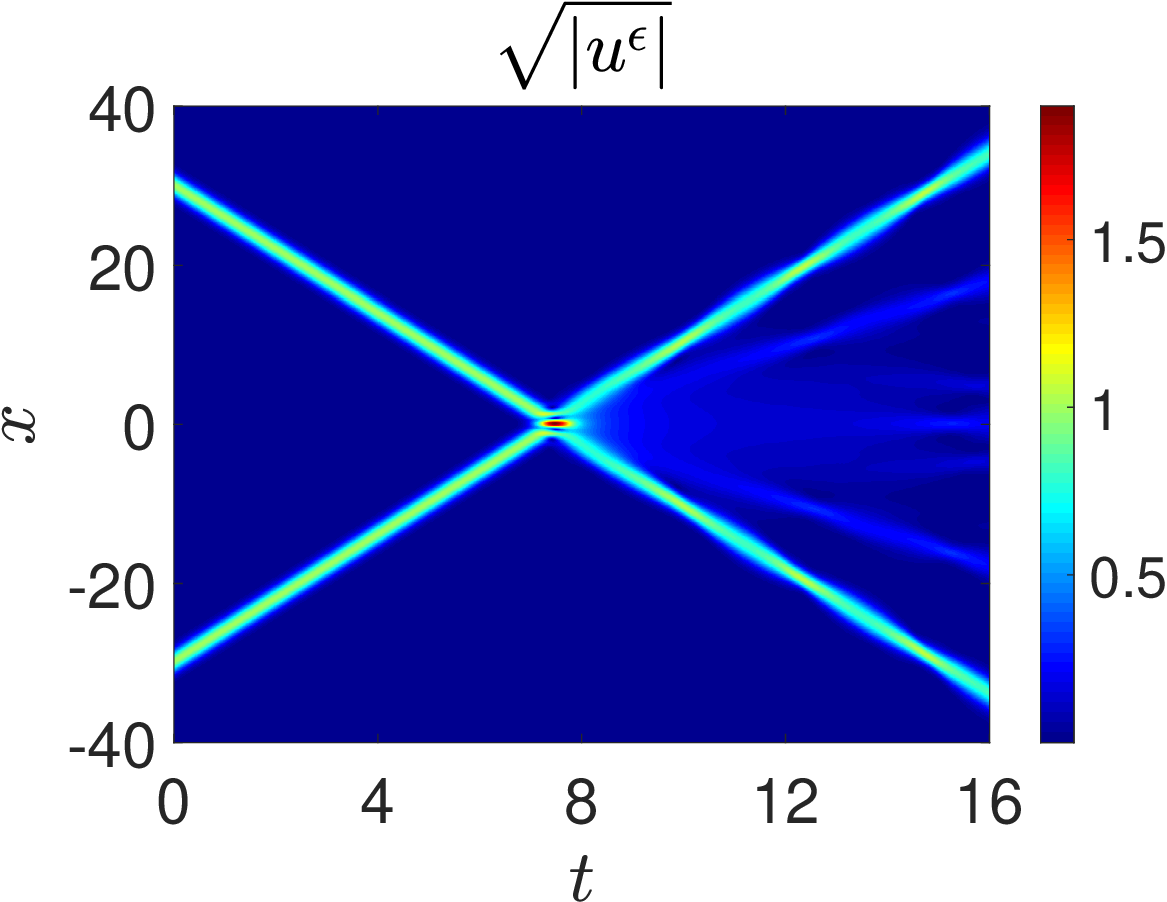}
\end{minipage}
\hspace{0.03\textwidth}
\begin{minipage}[t]{65mm}
\includegraphics[width=65mm]{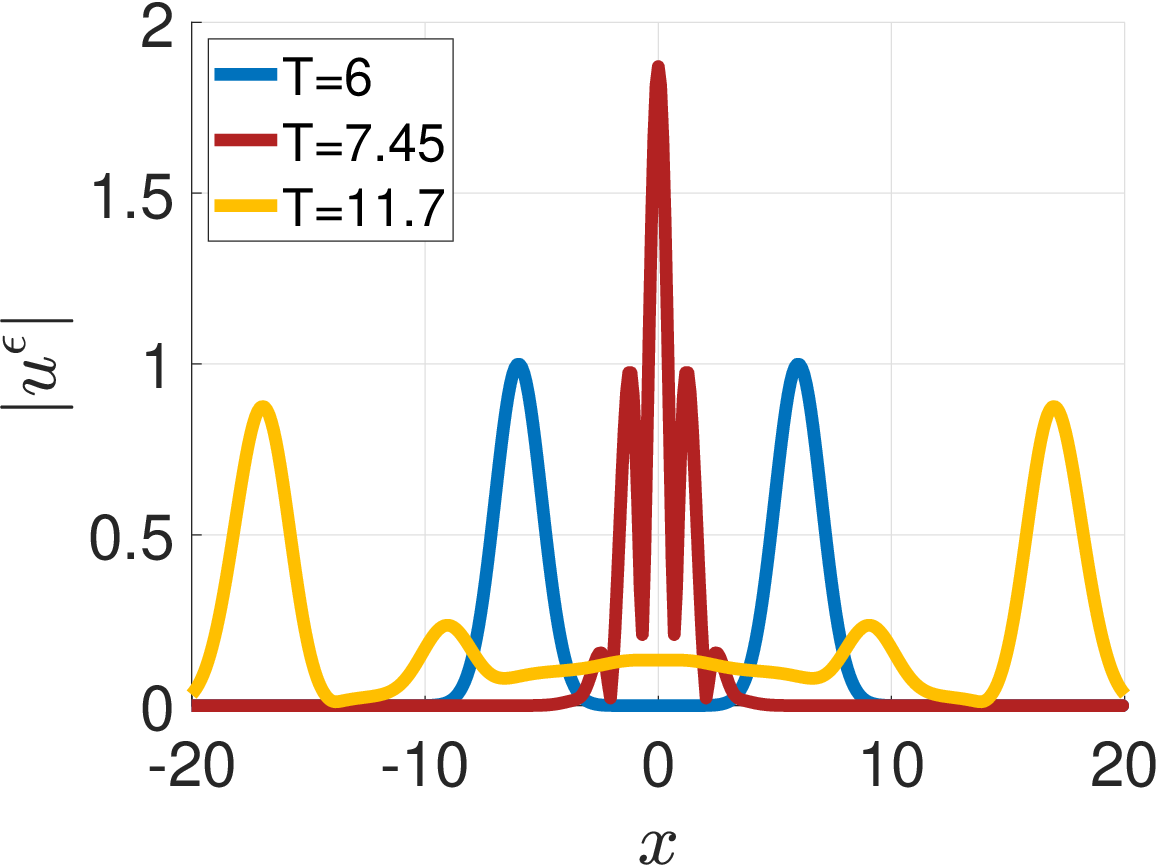}
\end{minipage}
\caption{Case III : Evolution of $\sqrt{\left| u^\epsilon \right|}$ from $t = 0$ to $t = 16$ (left), and plots of $\left| u^\epsilon \right|$ at different times (right) produced by {\bf SVM4} in example \ref{EX4-2}.}\label{1d-scheme:hotcase3}
\end{figure}

\begin{figure}[H]
\centering\begin{minipage}[t]{65mm}
\includegraphics[width=65mm]{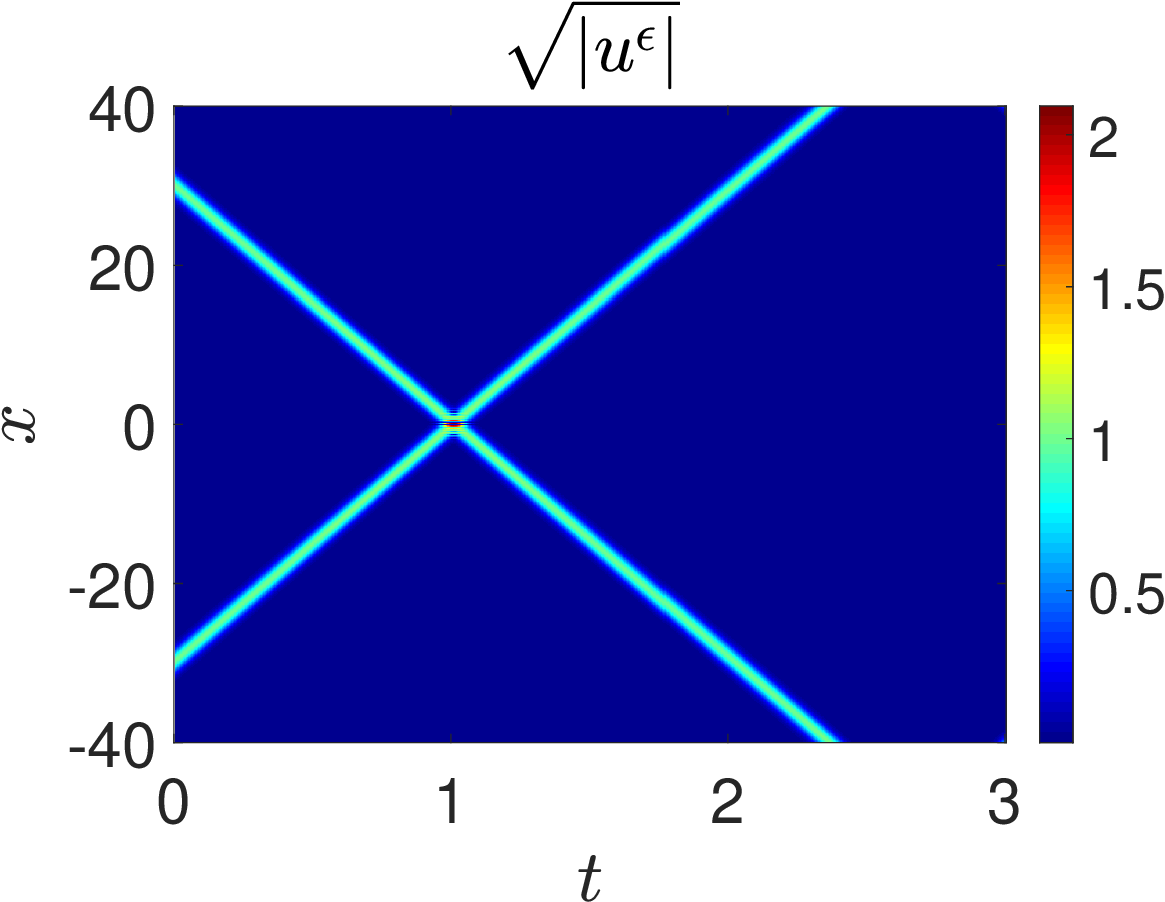}
\end{minipage}
\hspace{0.03\textwidth}
\begin{minipage}[t]{65mm}
\includegraphics[width=65mm]{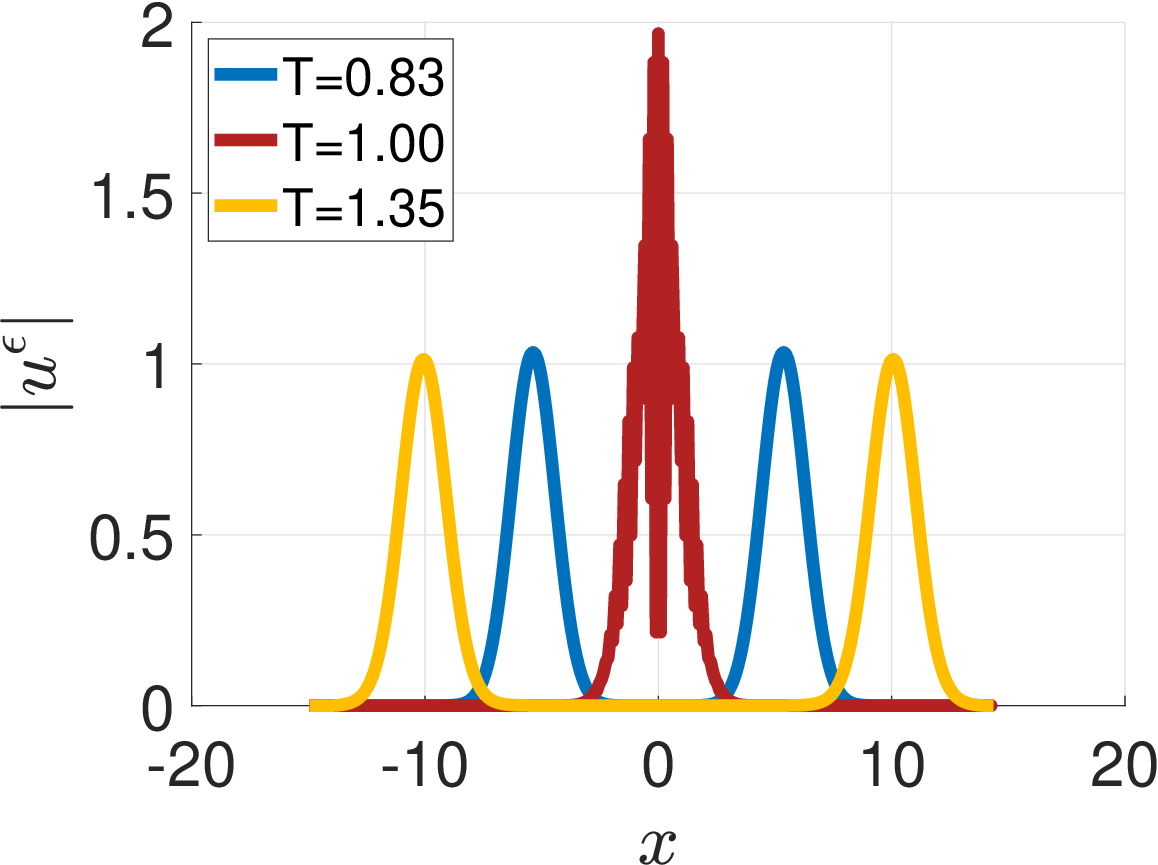}
\end{minipage}
\caption{Case IV : Evolution of $\sqrt{\left| u^\epsilon \right|}$ from $t = 0$ to $t = 3$ (left), and plots of $\left| u^\epsilon \right|$ at different times (right) produced by {\bf SVM4} in example \ref{EX4-2}.}\label{1d-scheme:hotcase4}
\end{figure}

Figures \ref{1d-scheme:hotcase1}-\ref{1d-scheme:hotcase4} display the evolution of Gaussons in 1D for Case
I-IV. In Figures \ref{1d-scheme:hotcase1} and \ref{1d-scheme:hotcase2}, it is clear to see that for initially well-separated static Gaussons, they will maintain unaltered  density profiles, while as two static Gaussians are close together, they will first move towards each other, collide and stick together, and then separate, swinging like a pendulum. Additionally, Figure \ref{1d-scheme:hotcase2} also shows that small outward isolated waves are emitted as the Gaussians separate and this pendulum motion of this emitted wave becomes faster as time goes on. When velocities are introduced to the Gaussons (i.e., $v_1=-v_2=2$), as illustrated in Figure \ref{1d-scheme:hotcase3}, we can observe that the two Gaussons undergo complete transmission through each other and ultimately move separately along with generating new Gaussons. However, if the velocities are increased from 2 (i.e., $v_1=-v_2=2$) to 15 (i.e., $v_1=-v_2=15$), Figure \ref{1d-scheme:hotcase4} demonstrates the two Gaussons continue to move towards each other at constant velocities, collide and separate accompanied by same move velocities and amplitudes. But an interesting  phoneme is that, different from Case III, no new Gaussons generate after collision. We note that the results agree well with those obtained by Qian et al.~\cite{qian2023}.
\begin{figure}[H]
\centering\begin{minipage}[t]{65mm}
\includegraphics[width=65mm]{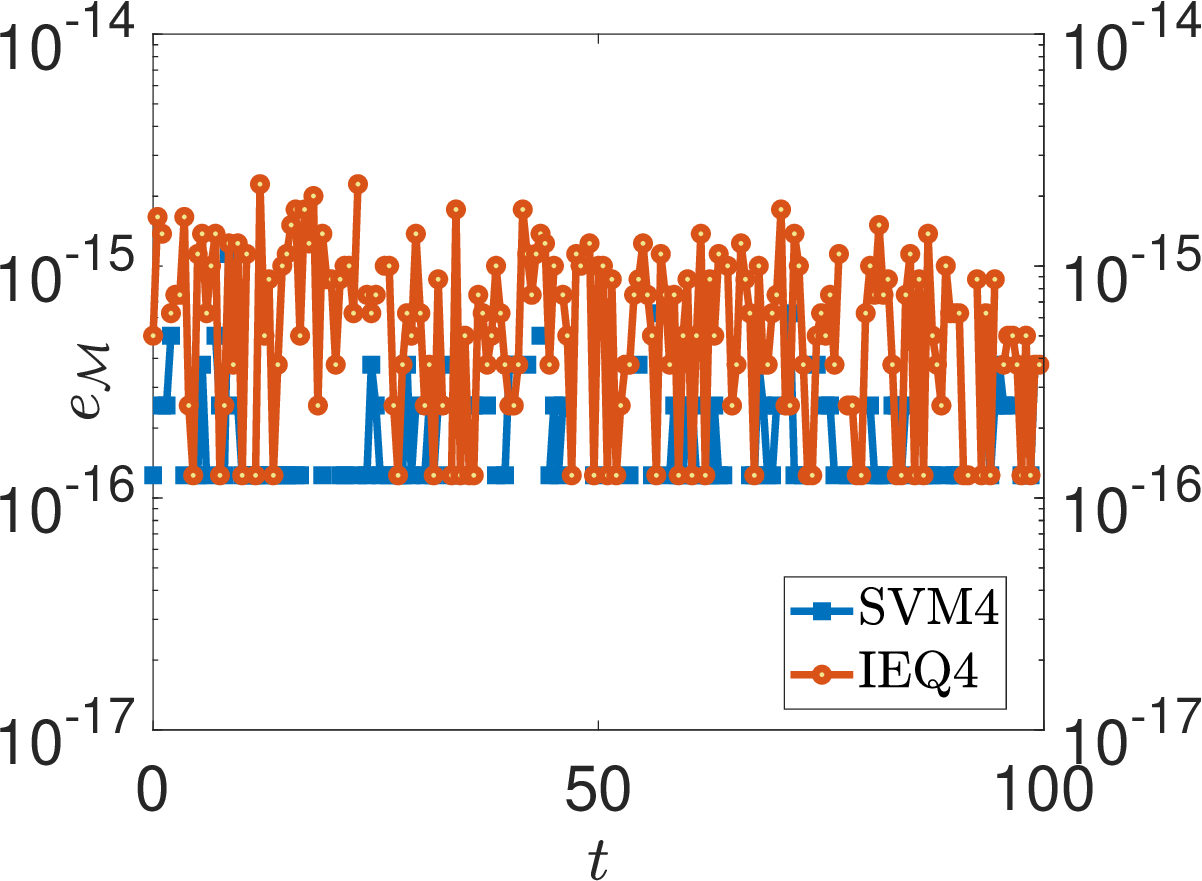}
\end{minipage}
\hspace{0.03\textwidth}
\begin{minipage}[t]{65mm}
\includegraphics[width=65mm]{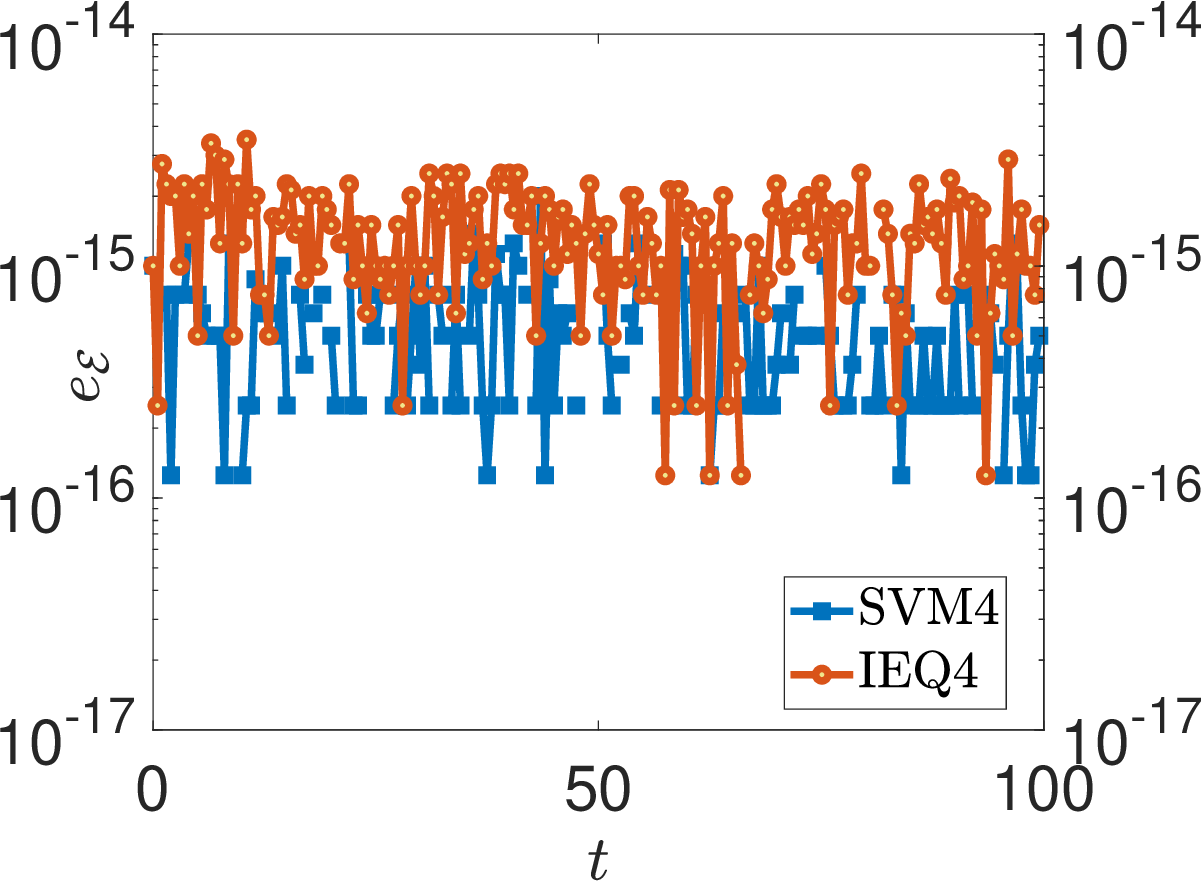}
\end{minipage}
\caption{Case I: The long-time evolution of $e_{\mathcal{M}}$ \& $e_{\mathcal{E}}$ for the two schemes in example \ref{EX4-2}.}\label{1d-scheme:fig:error_case1}
\end{figure}

\begin{figure}[H]
\centering\begin{minipage}[t]{65mm}
\includegraphics[width=65mm]{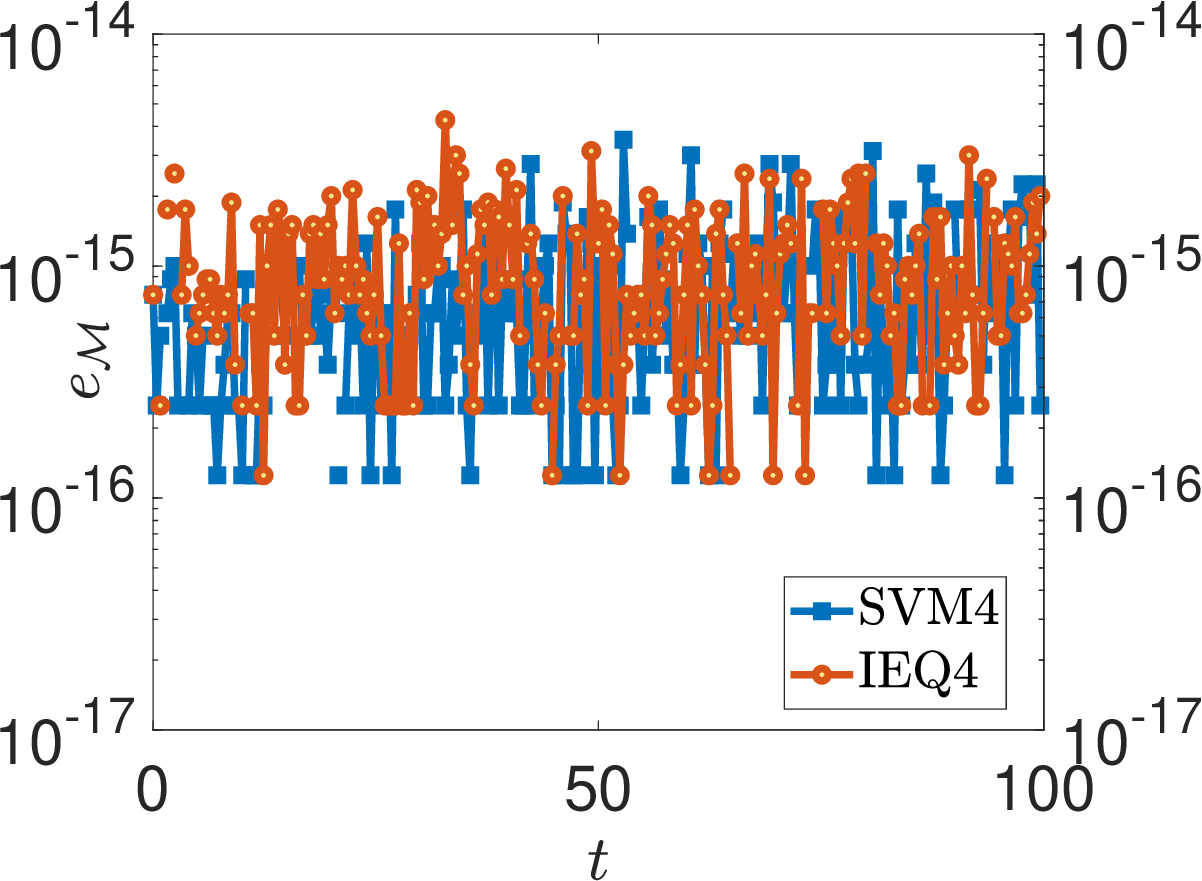}
\end{minipage}
\hspace{0.03\textwidth}
\begin{minipage}[t]{65mm}
\includegraphics[width=65mm]{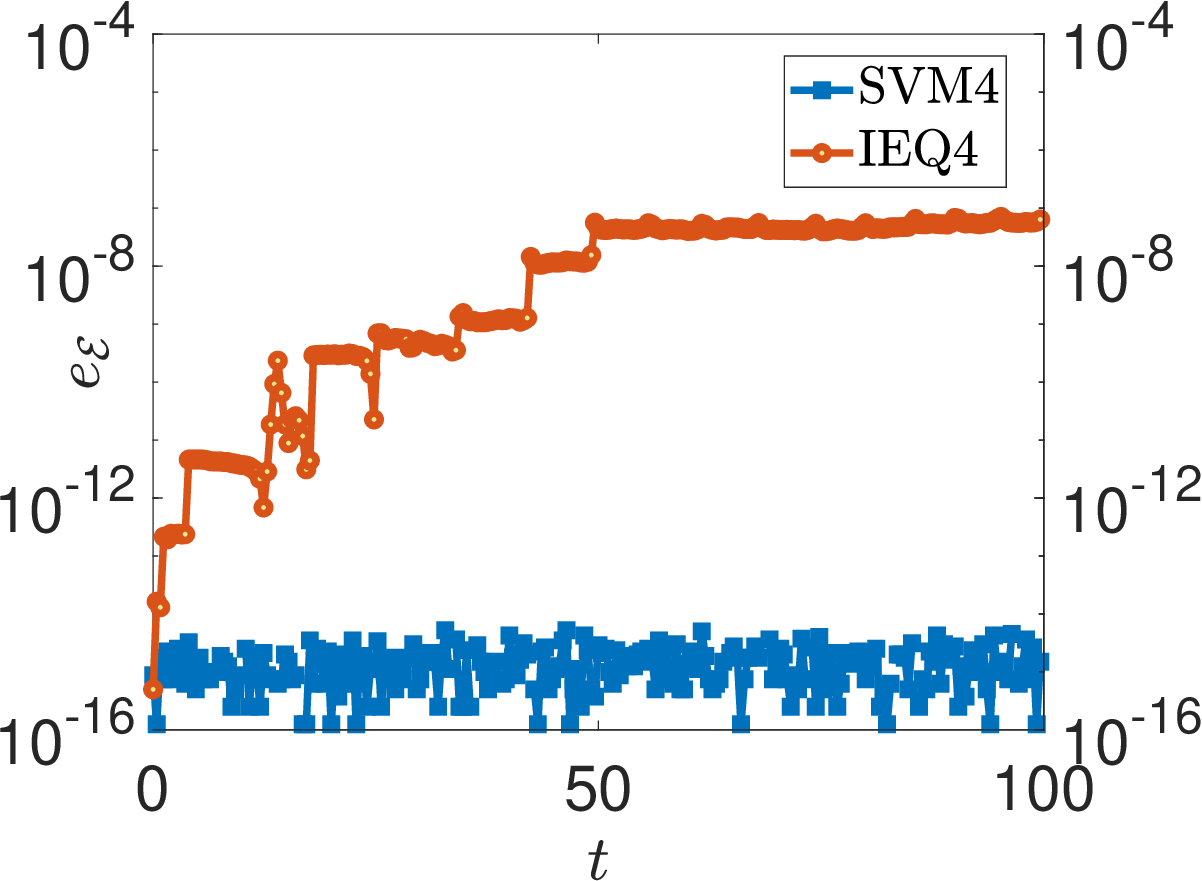}
\end{minipage}
\caption{Case II: The long-time evolution of $e_{\mathcal{M}}$ \& $e_{\mathcal{E}}$ for the two schemes in example \ref{EX4-2}.}\label{1d-scheme:fig:error_case2}
\end{figure}

\begin{figure}[H]
\centering\begin{minipage}[t]{65mm}
\includegraphics[width=65mm]{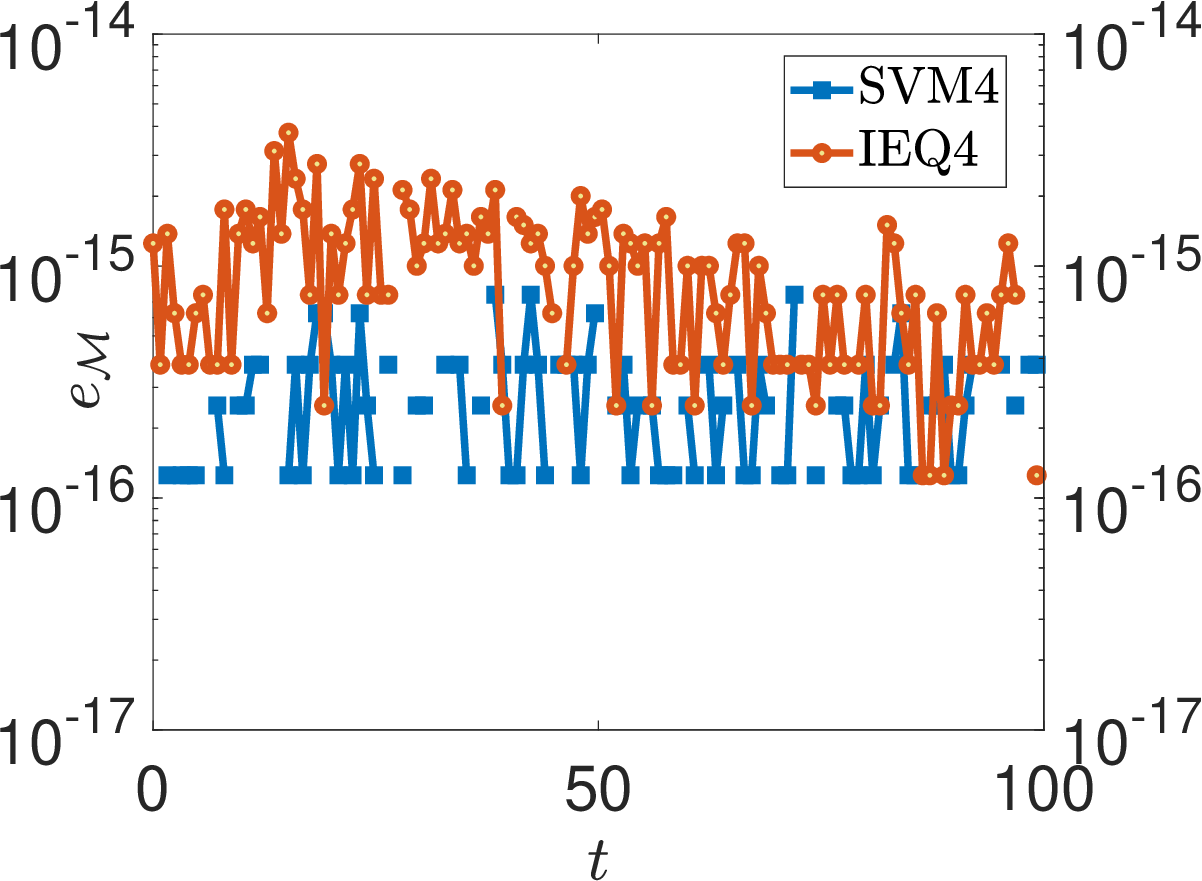}
\end{minipage}
\hspace{0.03\textwidth}
\begin{minipage}[t]{65mm}
\includegraphics[width=65mm]{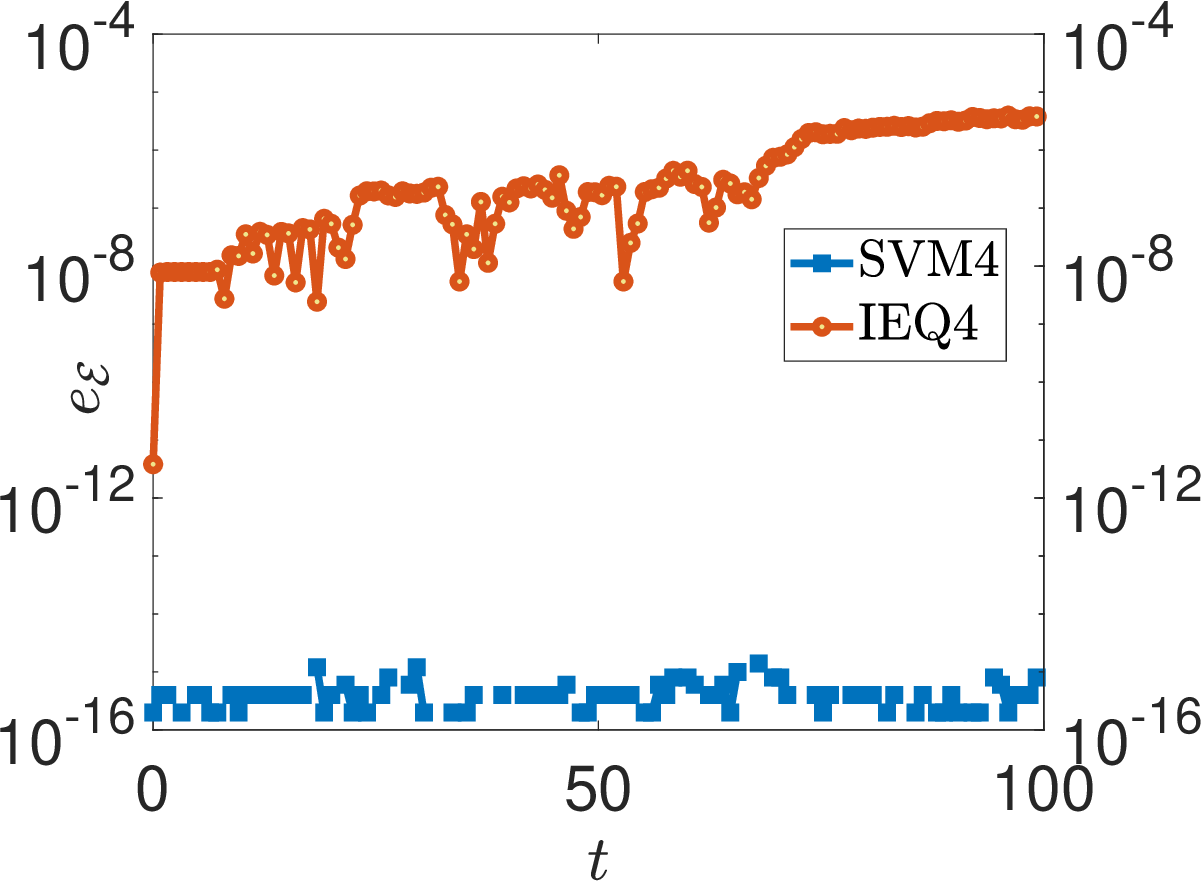}
\end{minipage}
\caption{Case III: The long-time evolution of $e_{\mathcal{M}}$ \& $e_{\mathcal{E}}$ for the two schemes in example \ref{EX4-2}.}\label{1d-scheme:fig:error_case3}
\end{figure}

\begin{figure}[H]
\centering\begin{minipage}[t]{65mm}
\includegraphics[width=65mm]{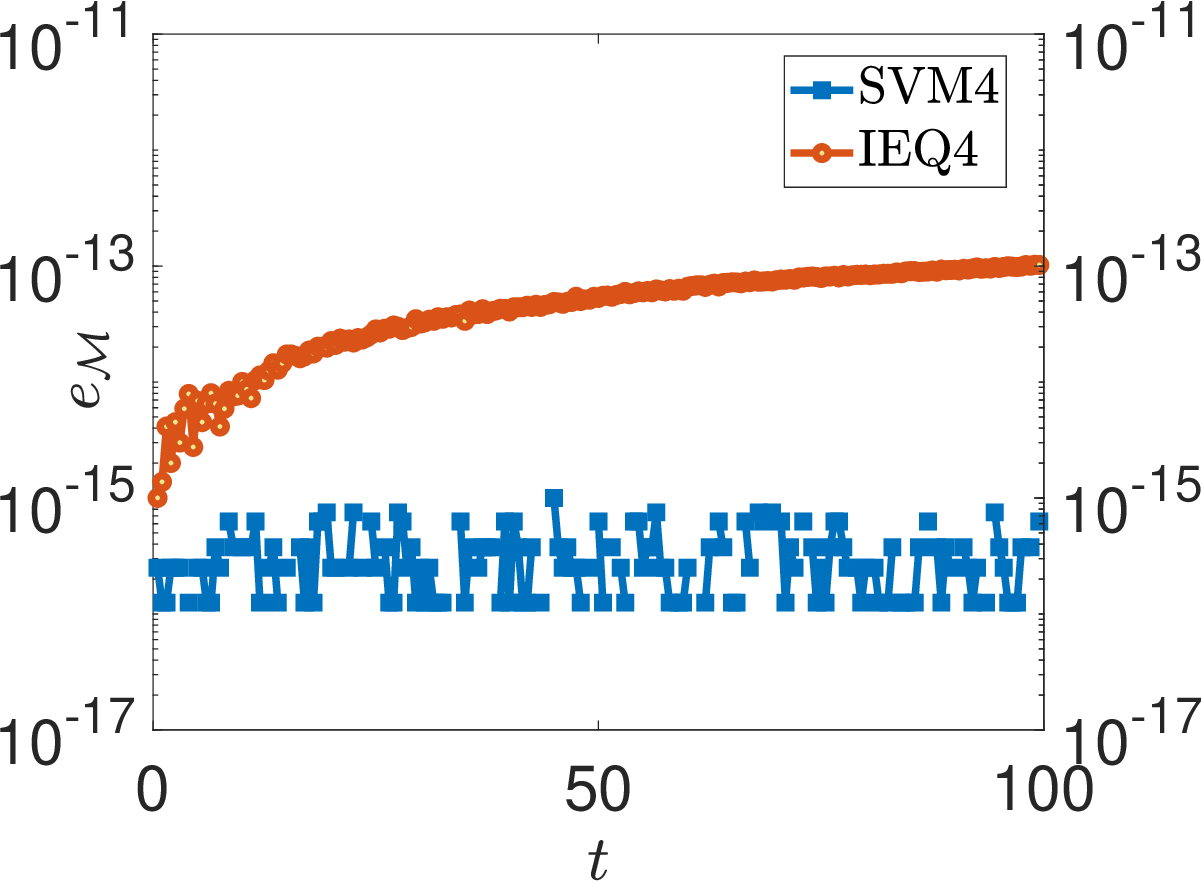}
\end{minipage}
\hspace{0.03\textwidth}
\begin{minipage}[t]{65mm}
\includegraphics[width=65mm]{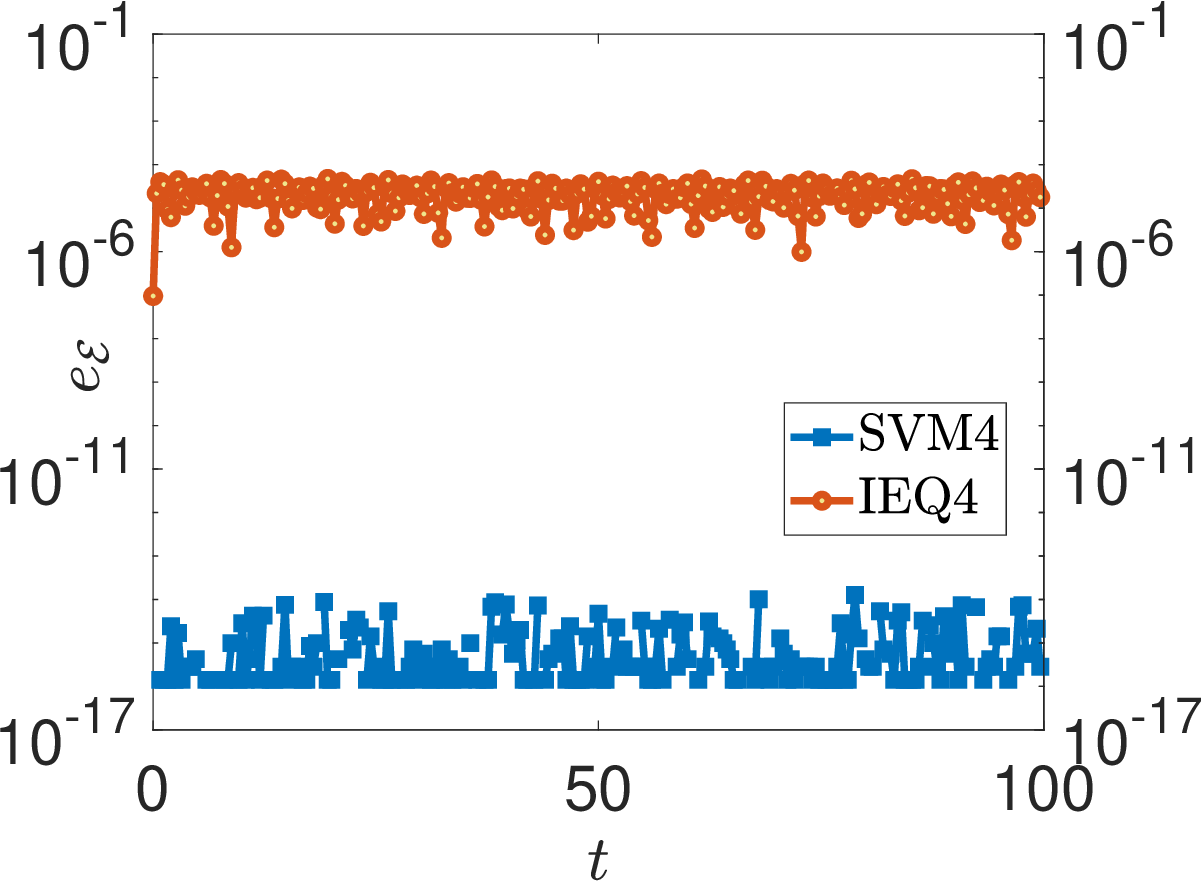}
\end{minipage}
\caption{Case IV: The long-time evolution of $e_{\mathcal{M}}$ \& $e_{\mathcal{E}}$ for the two schemes in example \ref{EX4-2}.}\label{1d-scheme:fig:error_case4}
\end{figure}
Figures \ref{1d-scheme:fig:error_case1}-\ref{1d-scheme:fig:error_case4} show the long-time evolution of $e_{\mathcal{M}}$ \& $e_{\mathcal{E}}$ of Case I-IV for the three schemes in example \ref{EX4-2} on the time interval $[0, 100]$. From Figures \ref{1d-scheme:fig:error_case1}-\ref{1d-scheme:fig:error_case4}, we can clearly observe that the errors on the mass and Hamiltonian energy produced by {\bf SVM4} are preserved to be machine precision, while the {\bf IEQ4} scheme only exactly preserve the mass conservation law.

\subsection{RlogSE in 2D}\label{sec:4.2}

\begin{ex}[\textbf{Accuracy confirmation in 2D}]\label{exmp2d-accuracy} In this example, we will test temporal numerical error, convergence order of the {\bf SVM4} and \textbf{IEQ4}scheme for the wave function ${u}^{\epsilon}(\cdot ,T=1)$ of the RlogSE \eqref{RLogSE} in 2D by taking $x_{L}=-x_{R}=-16,y_{L}=-y_{R}=-16$, the Fourier node $512^2$, the parameters $\lambda = -1,\ \epsilon = 1.0\ \times 10^{-12}$ and the following initial condition~\cite{BCST-MMMAS-2022}
\begin{align}
\label {2d-ex1}
u_{0}(\textbf{x})=\frac{1}{\sqrt [3]{\pi}}e^{\lambda/2|\textbf{x}-\textbf{x}^{0}|^{2}+\mathrm{i}\textbf{v}\textbf{x}},\ \textbf{x}\in\Omega,
\end{align}
where $\textbf{v}=(1, 1)^{\top}, \textbf{x}^{0}=(-2, 0)^{\top}$. In addition, we take the numerical solution produced by the proposed {\bf SVM4} with the time step $\tau = 2 \times 10^{-4}$  as a ``reference solution".
\end{ex}
 Table \ref{Tab:EX1-2D-order} reports the numerical errors and convergence orders. it is clearly demonstrated that the two schemes are fourth order accurate in time, and the errors produced by {\bf SVM4} are much smaller than the ones produced by {\bf IEQ4}. 

\begin{table}[H]
\tabcolsep=9pt
\footnotesize
\renewcommand\arraystretch{1.1}
\centering
\caption{Temporal error and order of convergence at $T = 1$ in example \ref{exmp2d-accuracy}. }
\label{Tab:EX1-2D-order}
\begin{tabular*}{\textwidth}[h]{@{\extracolsep{\fill}}c c c c c c c}\hline
{Scheme\ \ } &{}&{$\tau_0=1/40$} &{$\tau_0/2$} &{$\tau_0/4$} &{$\tau_0/8$} &{$\tau_0/16$} \\     
\hline
\multirow{2}{*}{SVM4}
  &{$e^{\epsilon,\tau,h}_{2}$}&{2.29e-04}&{1.44e-05}&{8.99e-07}&{5.62e-08}&{3.51e-09}
  \\[1ex]
 {}  &{Order}  &{-}& {3.996}&{3.999}&{4.000}&{4.000}\\[1ex]
\multirow{2}{*}{IEQ4}
   &{$e^{\epsilon,\tau,h}_{2}$}&{2.33e-04}&{1.46e-05}&{9.12e-07}&{5.70e-08}&{3.56e-09}
   \\[1ex]
{}  &{Order}  &{-}& {3.996}&{3.999}&{4.000}&{4.000}\\\hline
\end{tabular*}
\end{table}

\begin{ex}\label{EX4-4} [\textbf{Long-time evolution of dynamical behaviors and conservation laws in 2D}]  In this example, we will apply the {\bf SVM4} scheme to study the long time dynamical behaviors and conservation laws of the RlogSE \eqref{RLogSE} in 2D by taking the Fourier node $512^2$, the time step $\tau = 1 \times 10^{-2}$ with the parameters $\lambda = -1,\ \epsilon = 1.0\ \times 10^{-12}$ and the following initial condition~\cite{BCST-MMMAS-2022}
\begin{align}
\label {2d}
u_{0}(\textbf{x})=b_{1}e^{\lambda/2|\textbf{x}-\textbf{x}^{0}_{1}|^{2}+\mathrm{i}\textbf{v}_{1}\textbf{x}}+b_{2}e^{\lambda/2|\textbf{x}-
\textbf{x}^{0}_{2}|^{2}+\mathrm{i}\textbf{v}_{2}\textbf{x}},\ \textbf{x}\in\Omega,
\end{align}
where $ b_{k}, v_{k}$ and $\textbf{x}^{0}_{j},(j= 1, 2 )$  are real constant vectors. Here, we consider the following cases:
\begin{itemize}
\item Case I: $b_{1}=b_{2}=\frac{1}{\sqrt [4]{\pi}},\ \textbf{v}_{1}=\textbf{v}_{2}=(0, 0)^{\top}, \textbf{x}^{0}_{1}=-\textbf{x}^{0}_{2}=(-2, 0)^{\top}$;
\item Case II: $b_{1}=1.5 b_{2}=\frac{1}{\sqrt [4]{\pi}}, \textbf{v}_{1}=(-0.15, 0)^{\top}, \textbf{v}_{2}= \textbf{x}^{0}_{1}=(0, 0)^{\top}, \textbf{x}^{0}_{2}=(5, 0)^{\top}$;
\item Case III: $b_{1}=b_{2}=\frac{1}{\sqrt [4]{\pi}}, \textbf{v}_{1}=(0, 0)^{\top}, \textbf{v}_{2}= (0, 0.85)^{\top}, \textbf{x}^{0}_{1}=-\textbf{x}^{0}_{2}=(-2, 0)^{\top}$.
\end{itemize}
\end{ex}

\begin{figure}[H]
\centering
\begin{minipage}[t]{55mm}
\includegraphics[width=55mm]{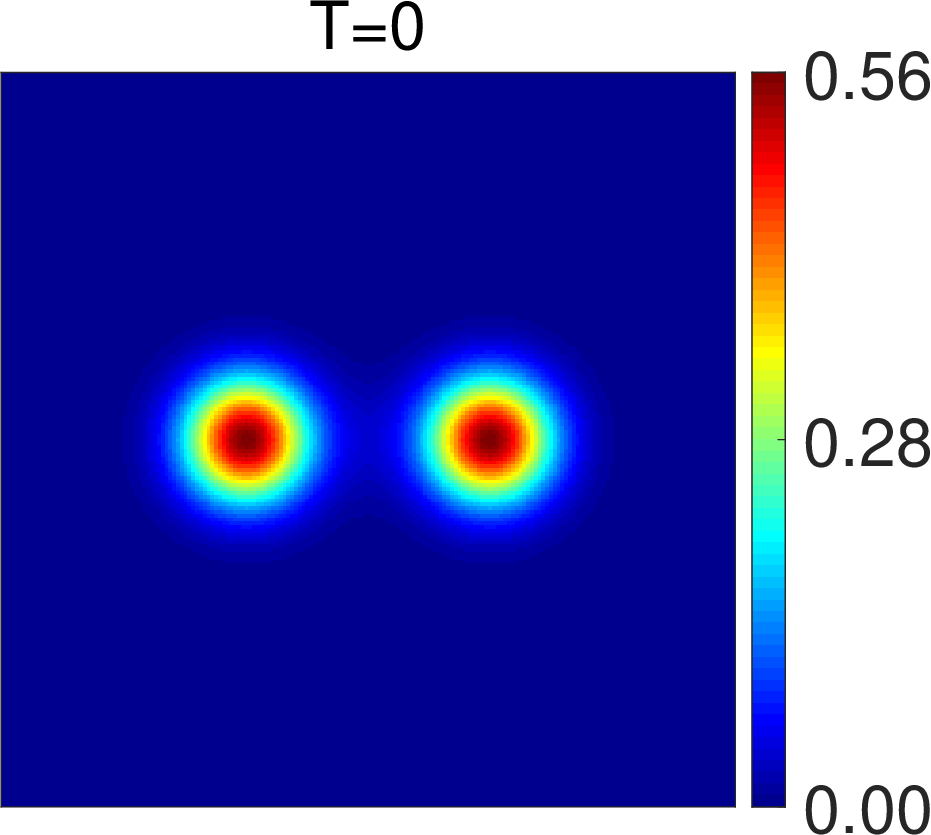}
\end{minipage}\
\hspace{0.03\textwidth}
\begin{minipage}[t]{55mm}
\includegraphics[width=55mm]{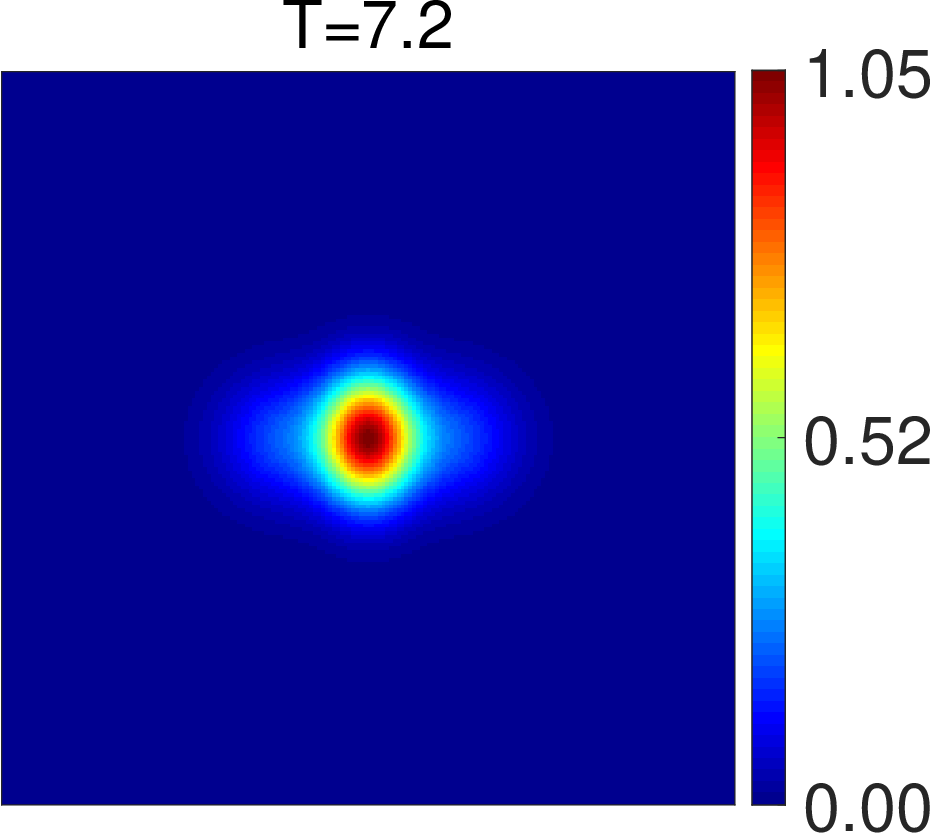}
\end{minipage}
\begin{minipage}[t]{55mm}
\includegraphics[width=55mm]{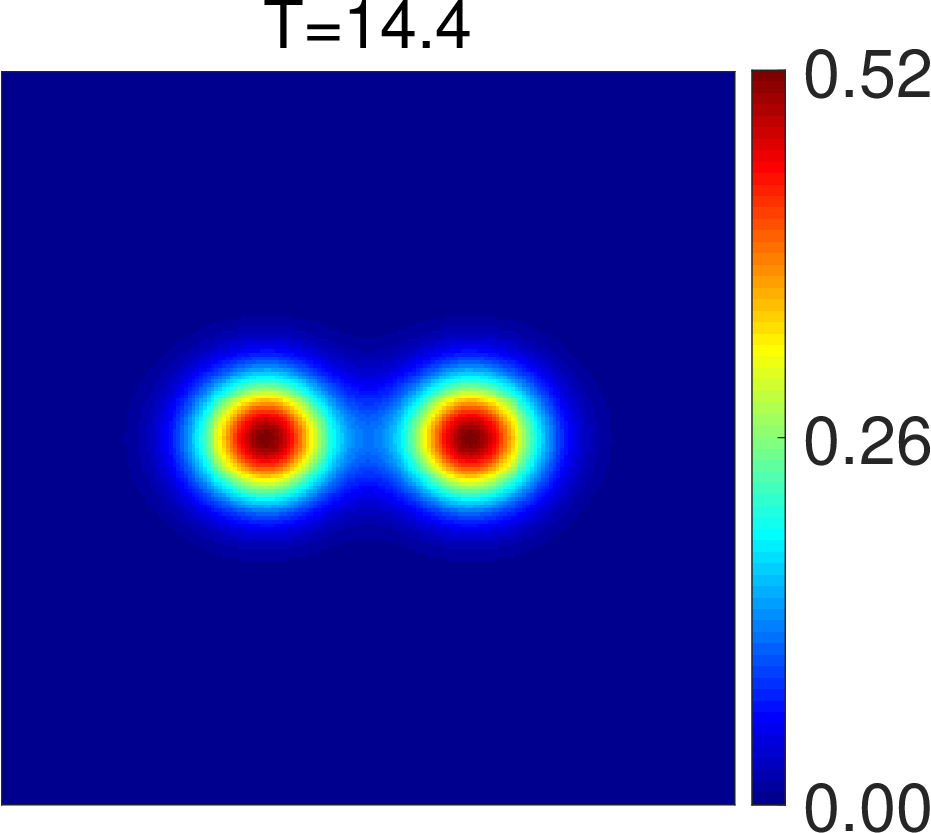}
\end{minipage}\
\hspace{0.03\textwidth}
\begin{minipage}[t]{55mm}
\includegraphics[width=55mm]{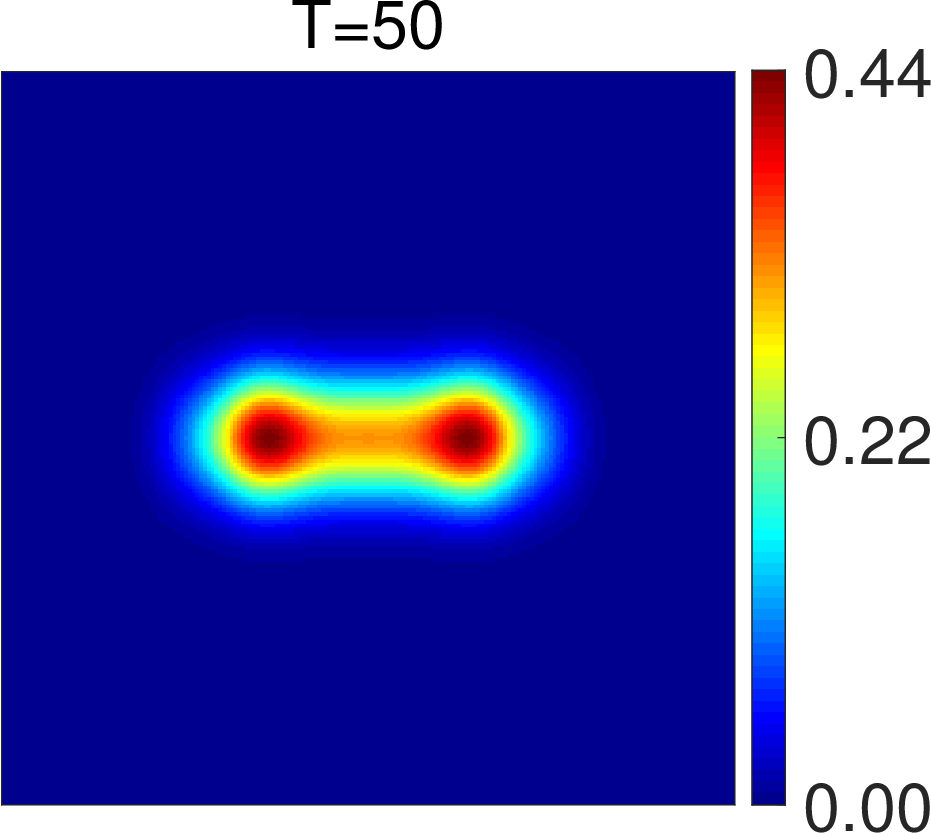}
\end{minipage}

\caption{Case I : Plots of $|u^{\epsilon}(\textbf{x},t)|^{2}$ produced by {\bf SVM4} at different times in spatial region $[-6, 6]^{2}$ in example \ref{EX4-4}.}\label{2d-scheme:fig:11}
\end{figure}

\begin{figure}[H]
\centering
\begin{minipage}[t]{55mm}
\includegraphics[width=55mm]{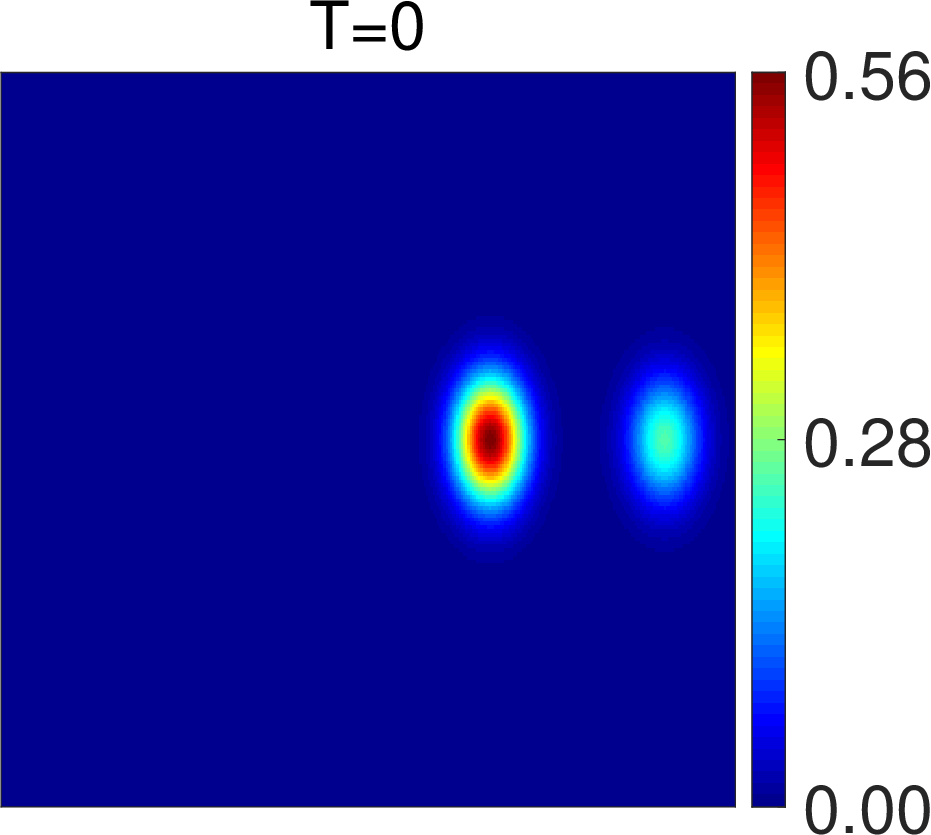}
\end{minipage}%
\hspace{0.03\textwidth}
\begin{minipage}[t]{55mm}
\includegraphics[width=55mm]{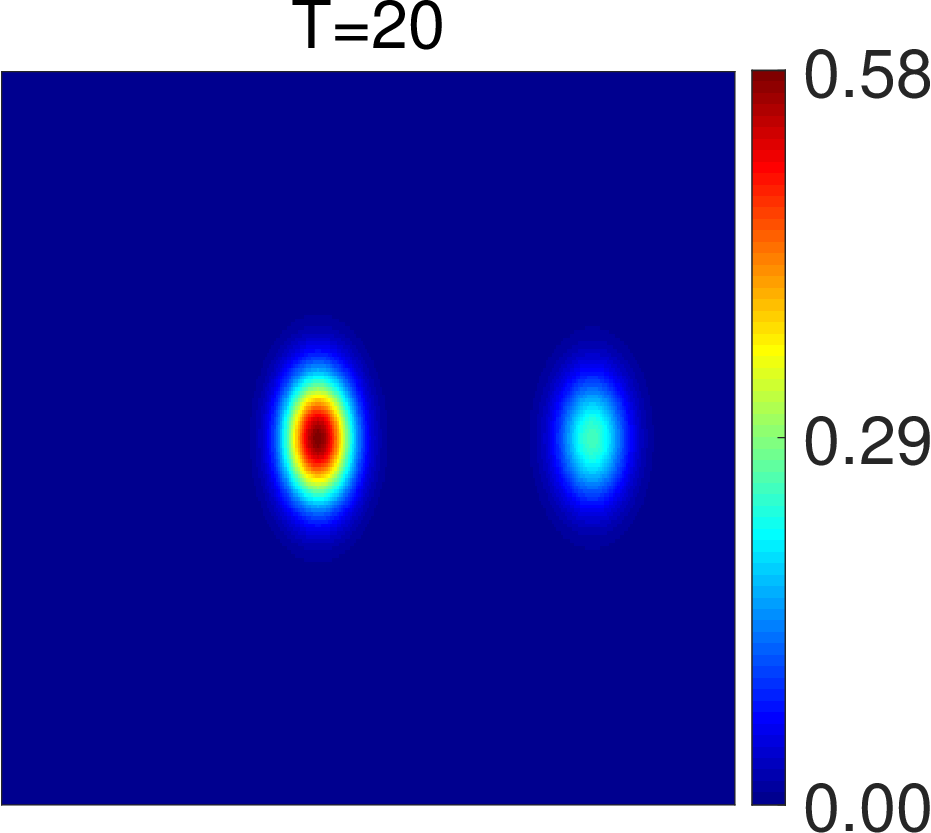}
\end{minipage}
\begin{minipage}[t]{55mm}
\includegraphics[width=55mm]{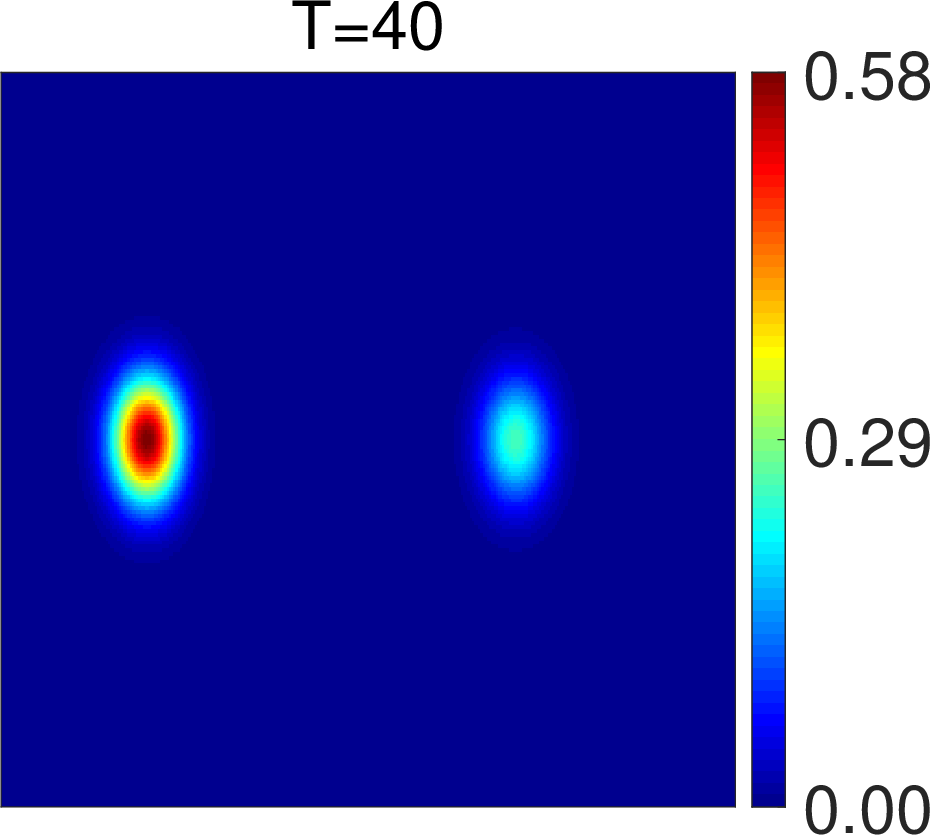}
\end{minipage}%
\hspace{0.03\textwidth}
\begin{minipage}[t]{55mm}
\includegraphics[width=55mm]{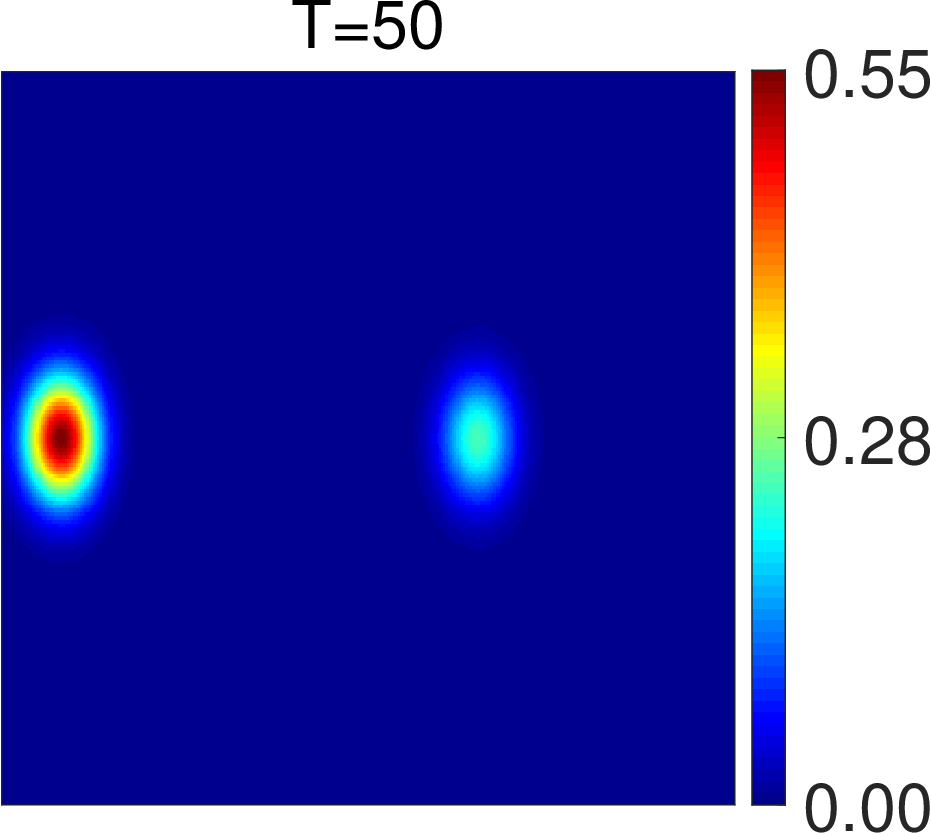}
\end{minipage}

\caption{Case II : Plots of $|u^{\epsilon}(\textbf{x},t)|^{2}$ produced by {\bf SVM4} at different times in spatial region $[-14, 7] \times [-6, 6]$ in example \ref{EX4-4}.}\label{2d-scheme:fig:13}
\end{figure}

\begin{figure}[H]
\centering
\begin{minipage}[t]{45mm}
\includegraphics[width=45mm]{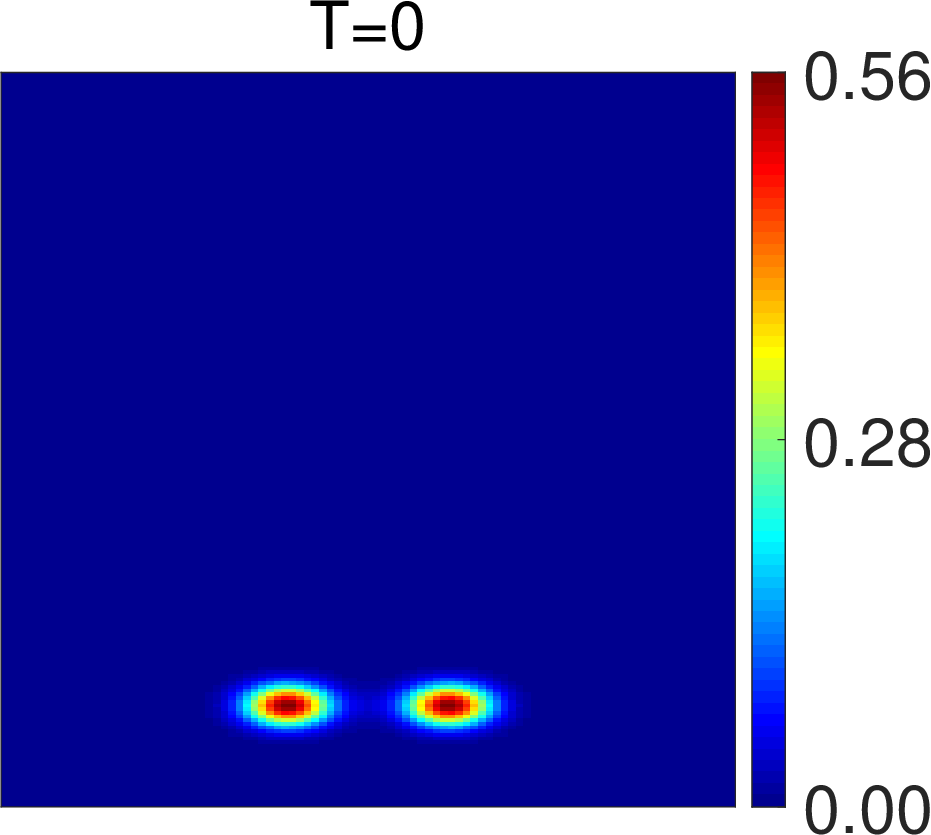}
\end{minipage}%
\hspace{0.02\textwidth}
\begin{minipage}[t]{45mm}
\includegraphics[width=45mm]{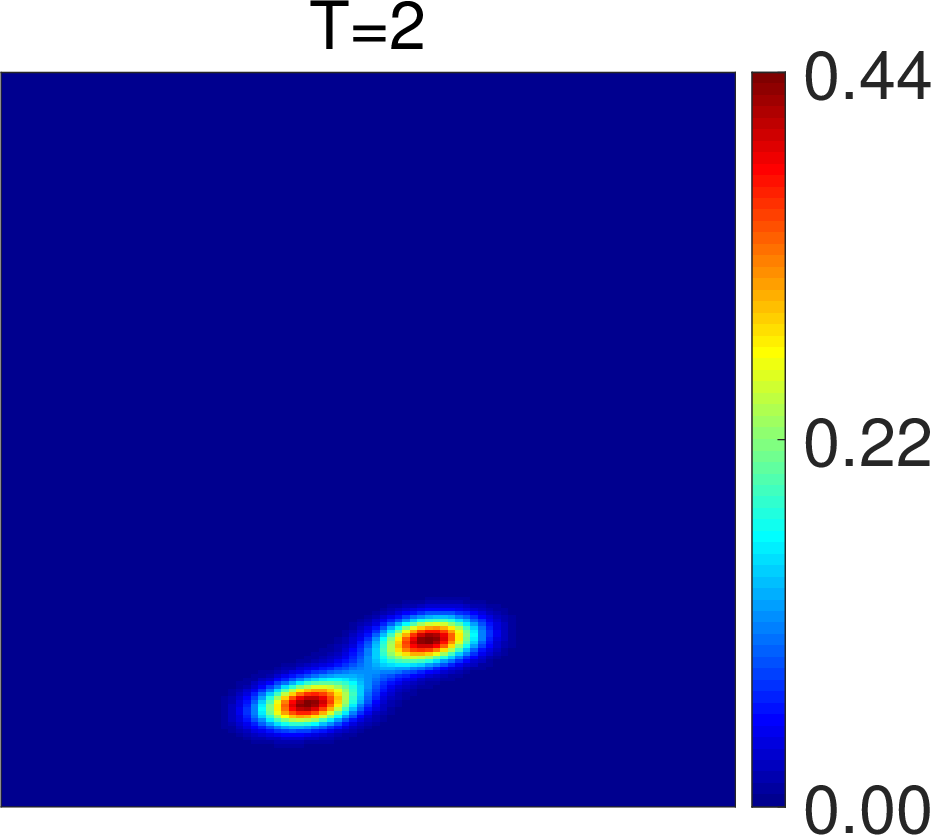}
\end{minipage}
\hspace{0.02\textwidth}
\begin{minipage}[t]{45mm}
\includegraphics[width=45mm]{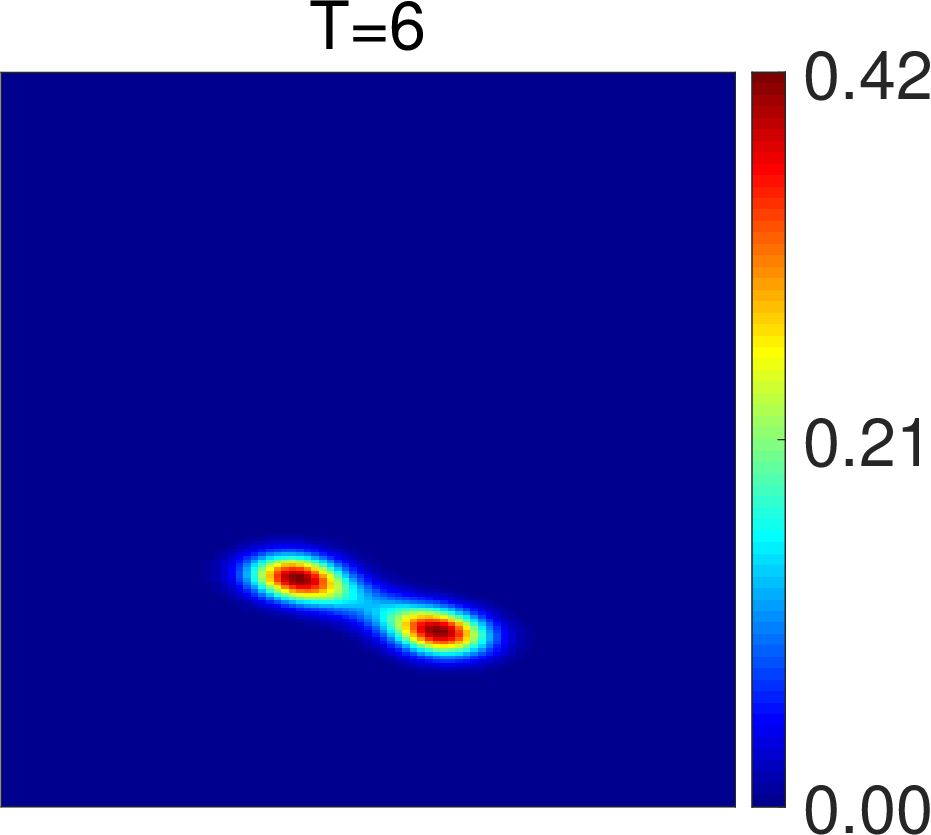}
\end{minipage}
\centering
\begin{minipage}[t]{45mm}
\includegraphics[width=45mm]{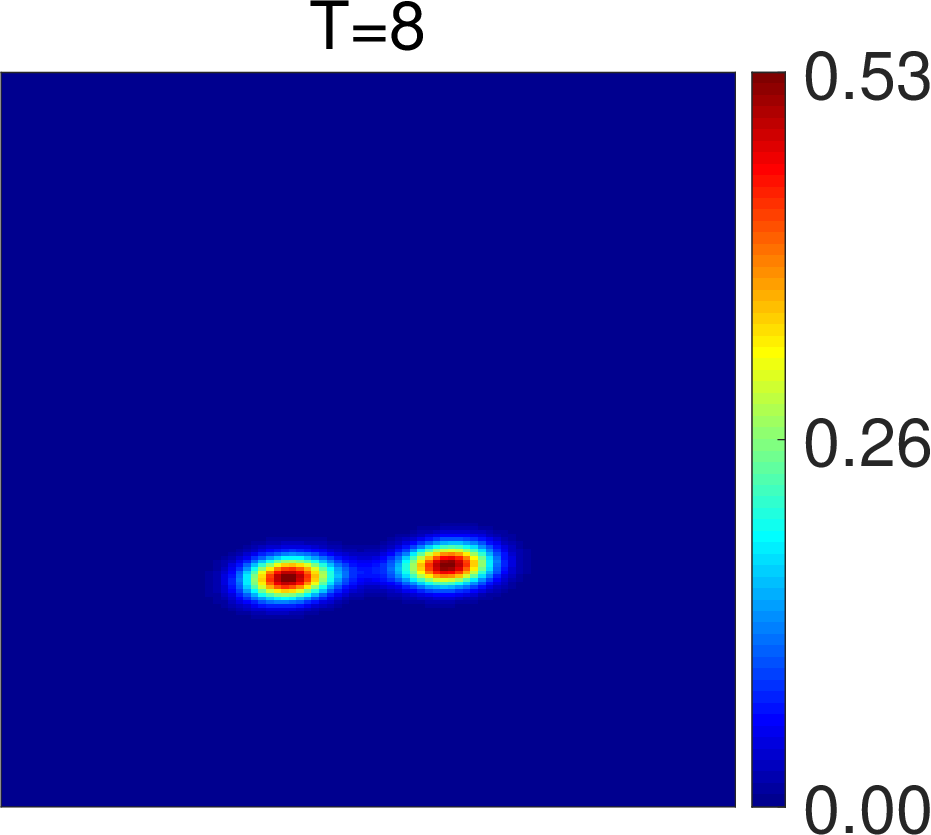}
\end{minipage}%
\hspace{0.02\textwidth}
\begin{minipage}[t]{45mm}
\includegraphics[width=45mm]{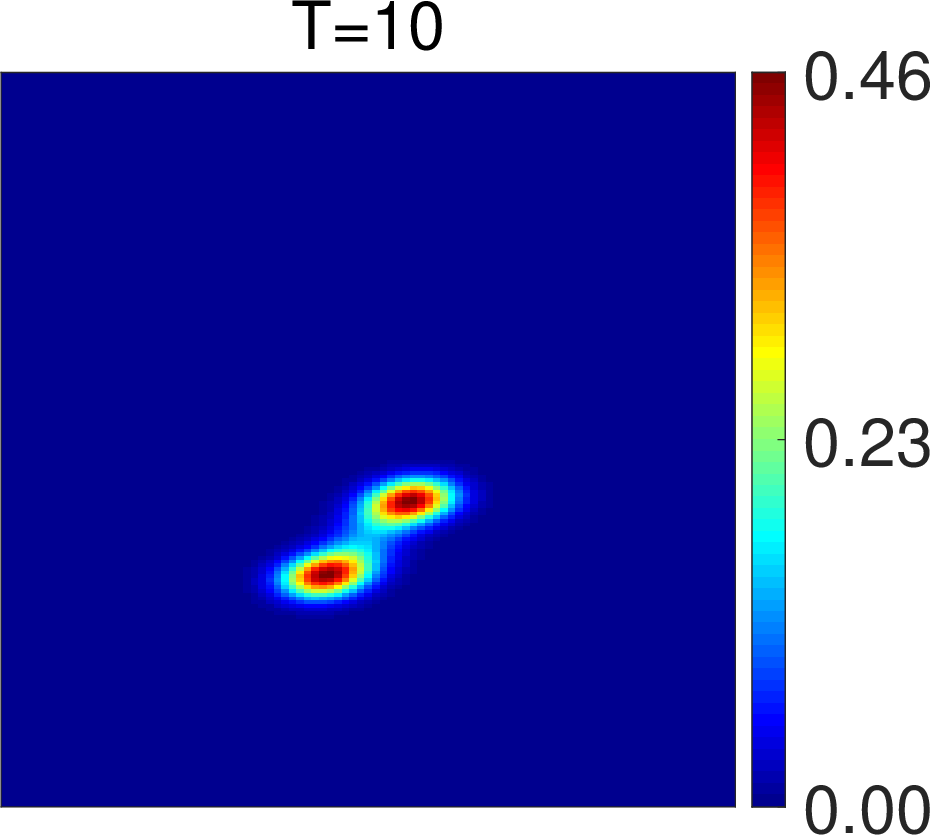}
\end{minipage}
\hspace{0.02\textwidth}
\begin{minipage}[t]{45mm}
\includegraphics[width=45mm]{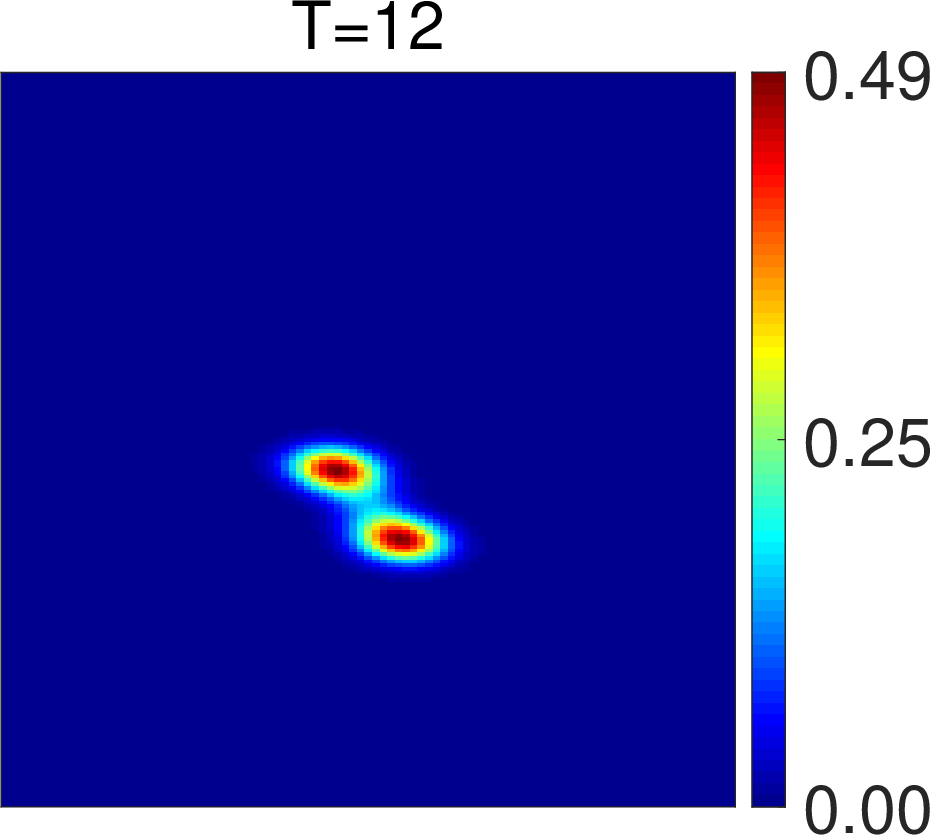}
\end{minipage}
\centering
\begin{minipage}[t]{45mm}
\includegraphics[width=45mm]{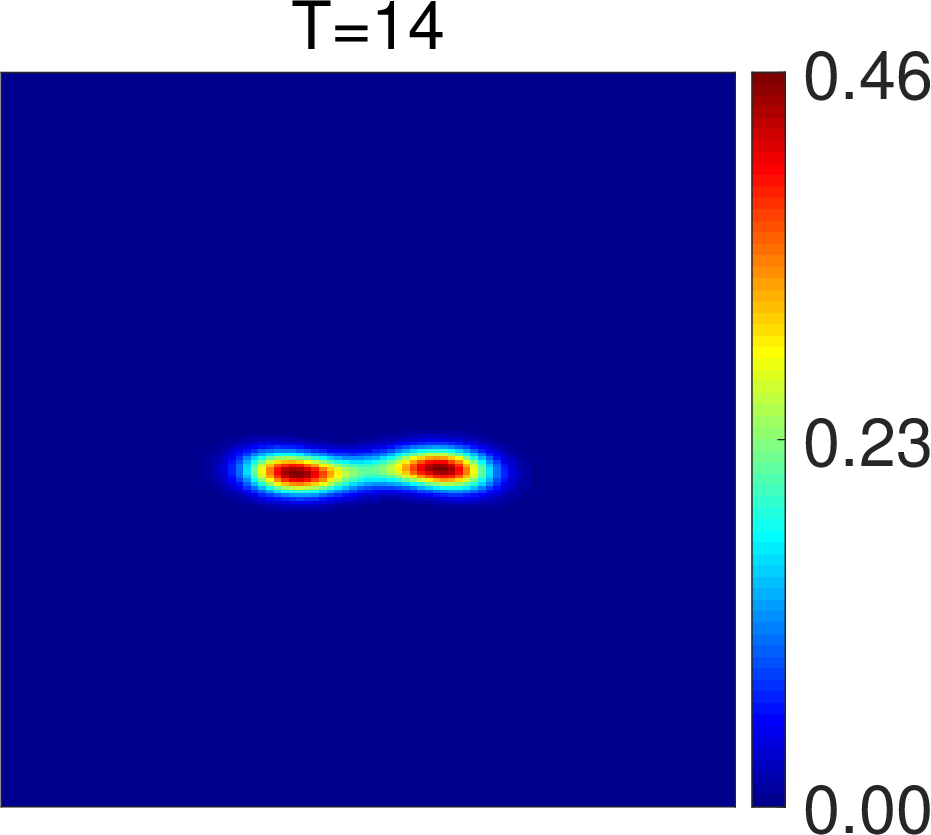}
\end{minipage}%
\hspace{0.02\textwidth}
\begin{minipage}[t]{45mm}
\includegraphics[width=45mm]{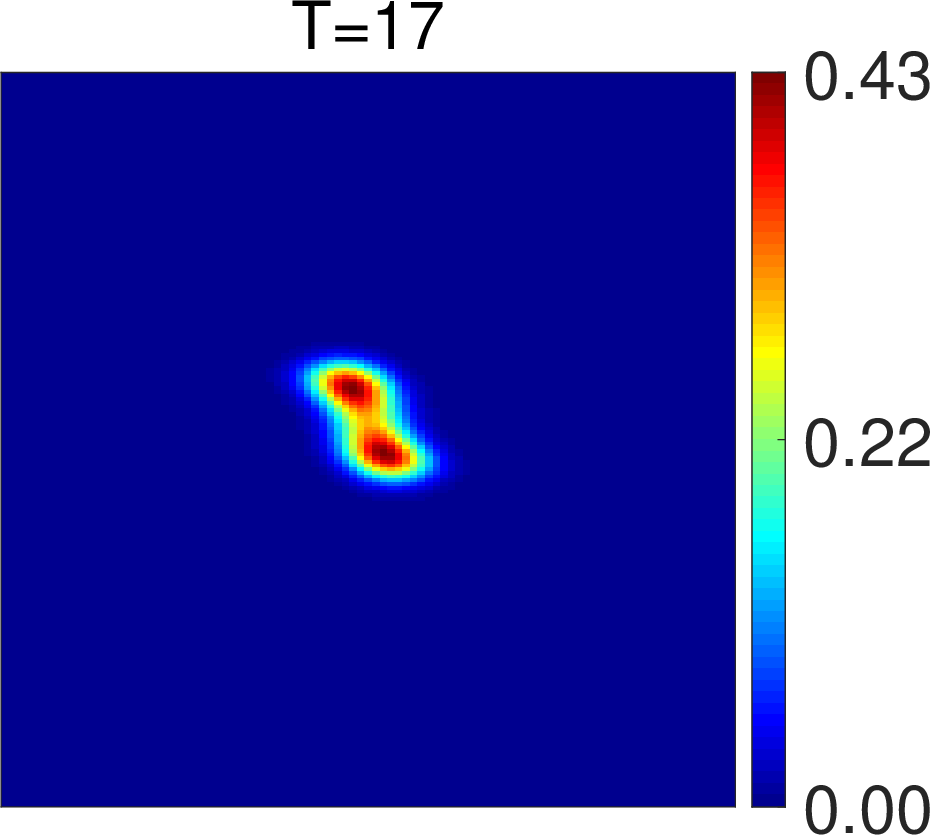}
\end{minipage}
\hspace{0.02\textwidth}
\begin{minipage}[t]{45mm}
\includegraphics[width=45mm]{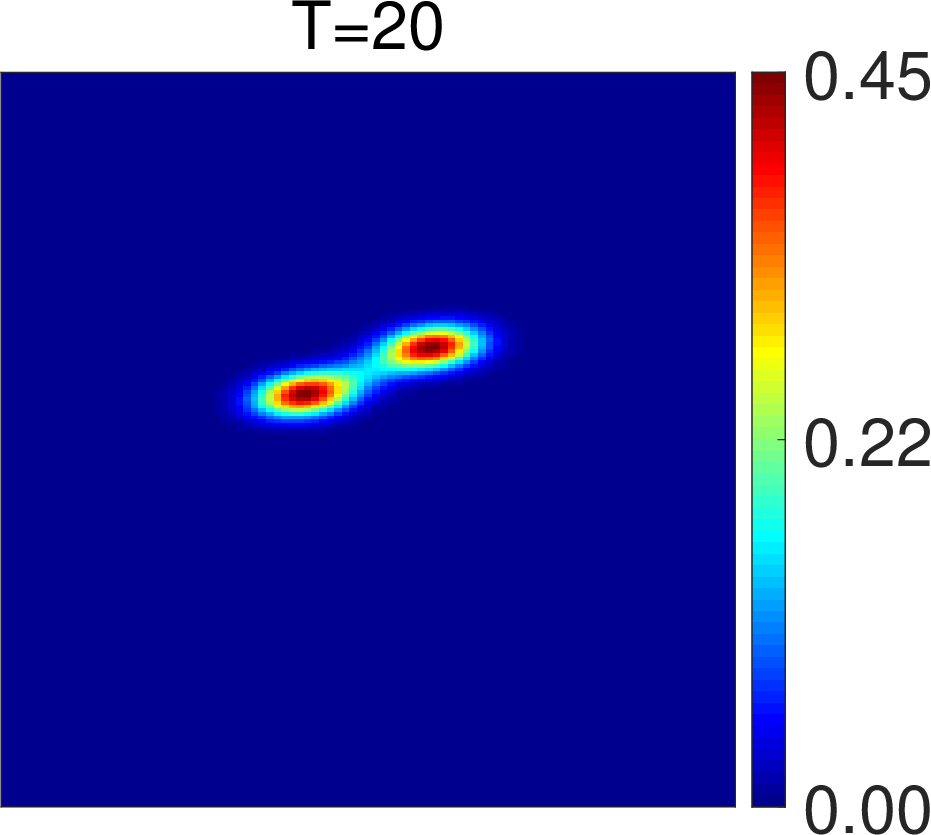}
\end{minipage}
\centering
\begin{minipage}[t]{45mm}
\includegraphics[width=45mm]{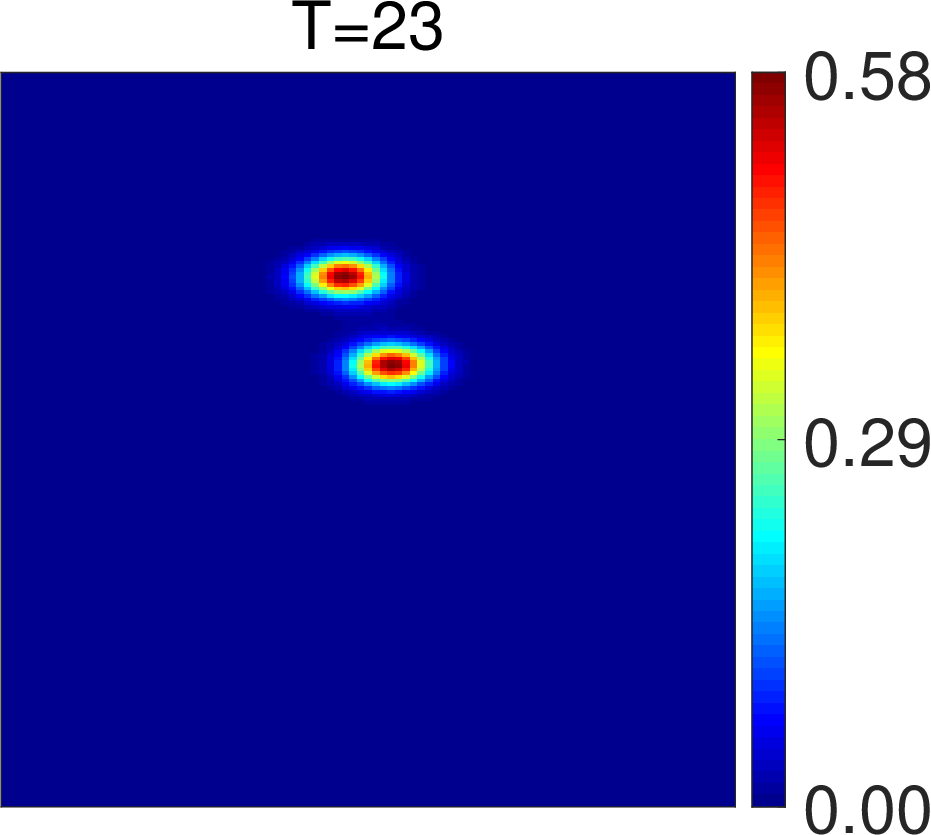}
\end{minipage}%
\hspace{0.02\textwidth}
\begin{minipage}[t]{45mm}
\includegraphics[width=45mm]{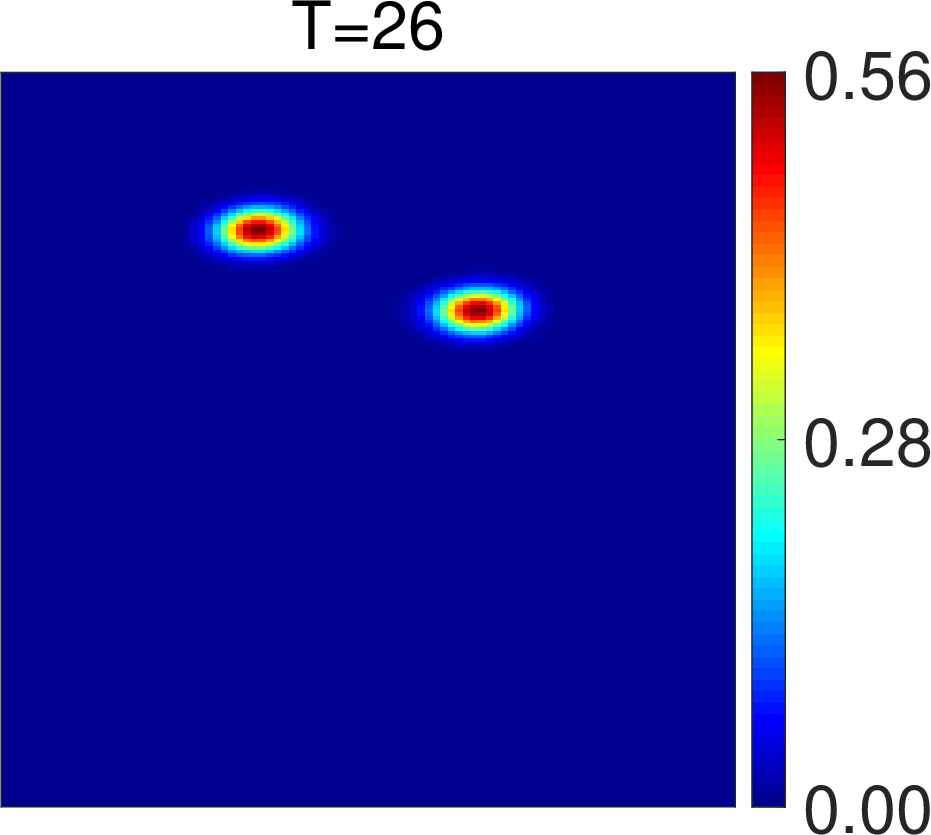}
\end{minipage}
\hspace{0.02\textwidth}
\begin{minipage}[t]{45mm}
\includegraphics[width=45mm]{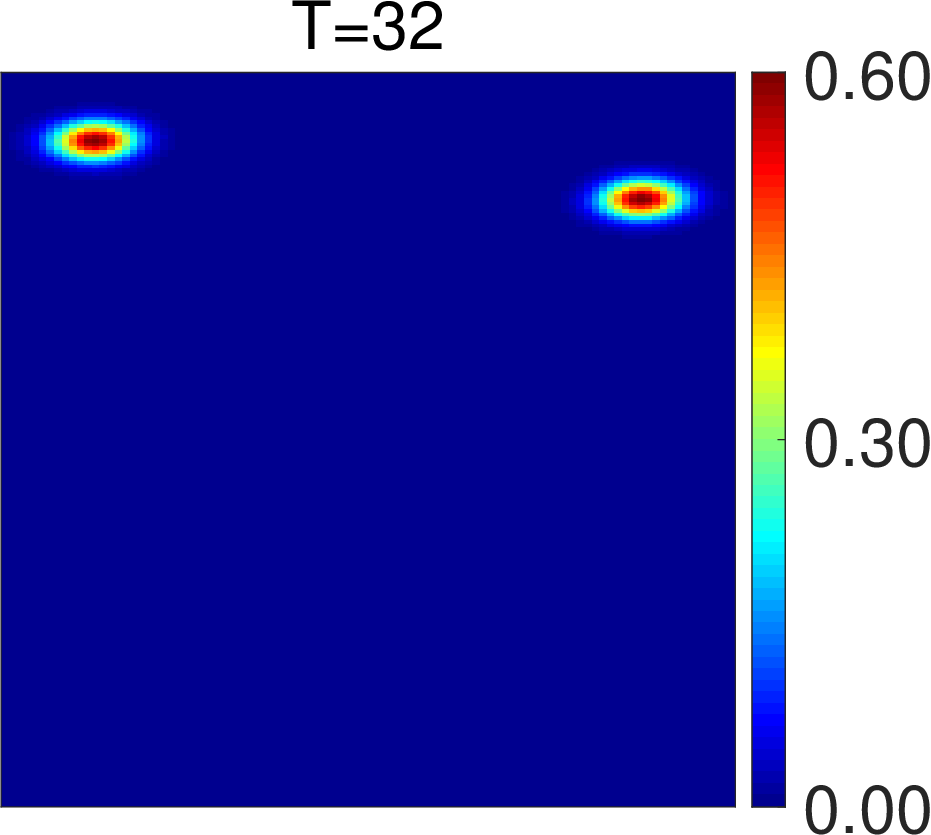}
\end{minipage}
\caption{Case III: Plots of $|u^{\epsilon}(\textbf{x},t)|^{2}$ produced by {\bf SVM4} at different times in spatial region $[-9, 9] \times [-5, 32]$ in example \ref{EX4-4}.}\label{2d-scheme:fig:15}
\centering\begin{minipage}[t]{65mm}
\includegraphics[width=65mm]{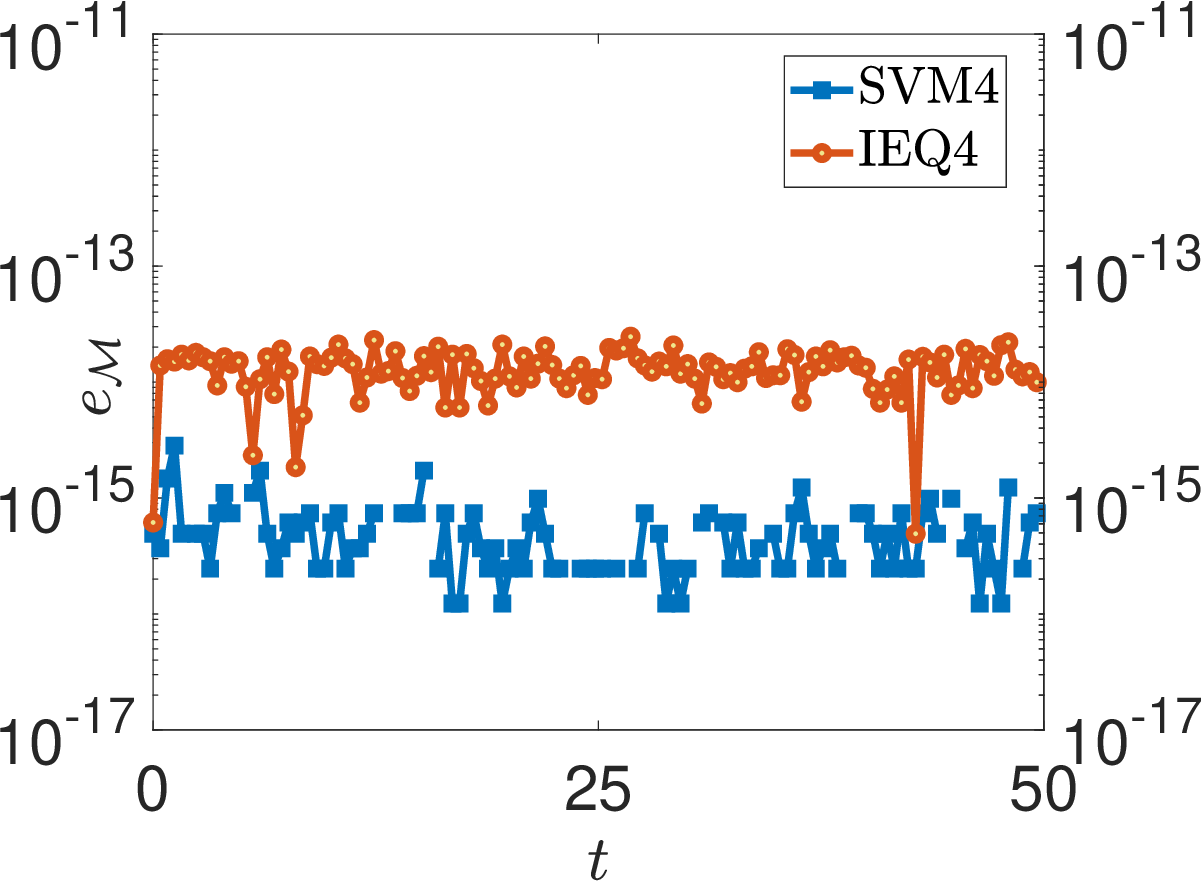}
\end{minipage}
\hspace{0.03\textwidth}
\begin{minipage}[t]{65mm}
\includegraphics[width=65mm]{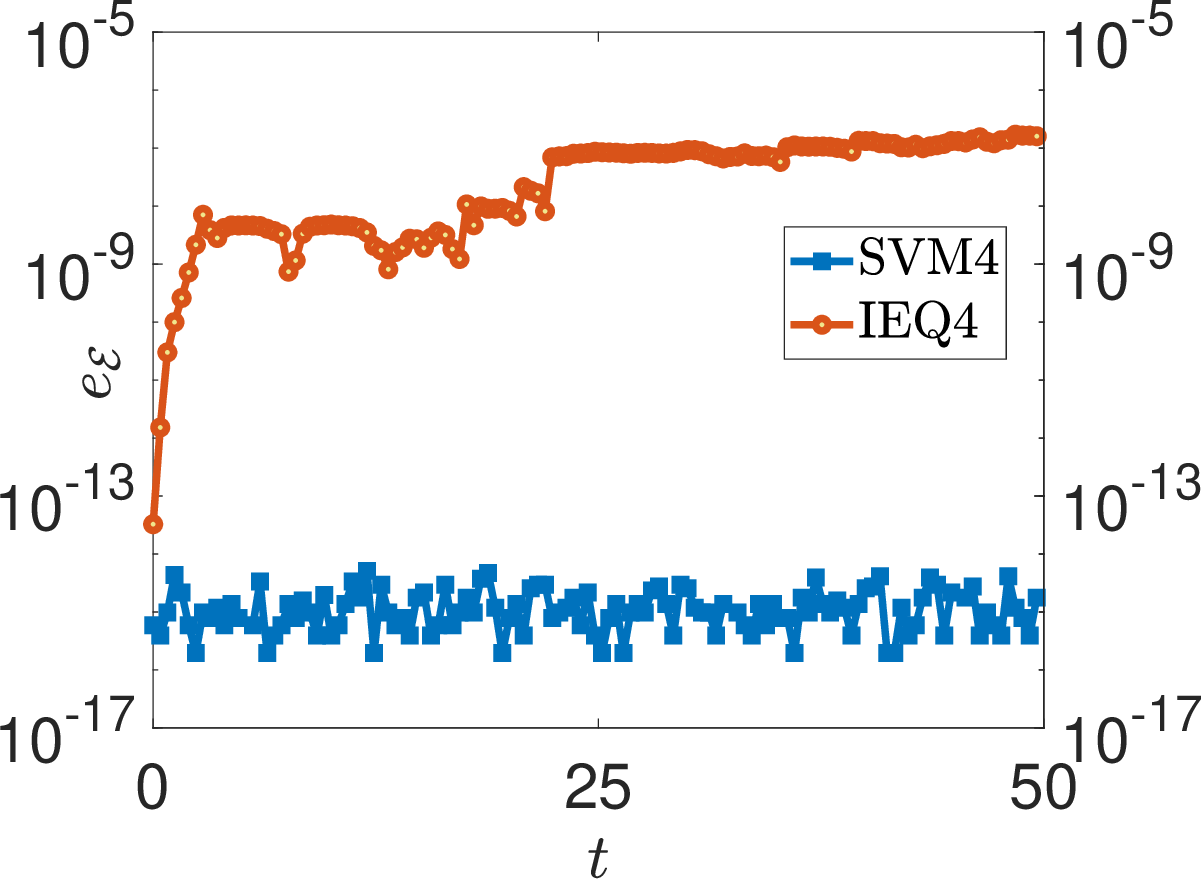}
\end{minipage}
\caption{Case I: The long-time evolution of $e_{\mathcal{M}}$ \& $e_{\mathcal{E}}$ for the two schemes in example \ref{EX4-4}.}\label{2d-scheme:case1}
\end{figure}

\begin{figure}[H]
\centering\begin{minipage}[t]{65mm}
\includegraphics[width=65mm]{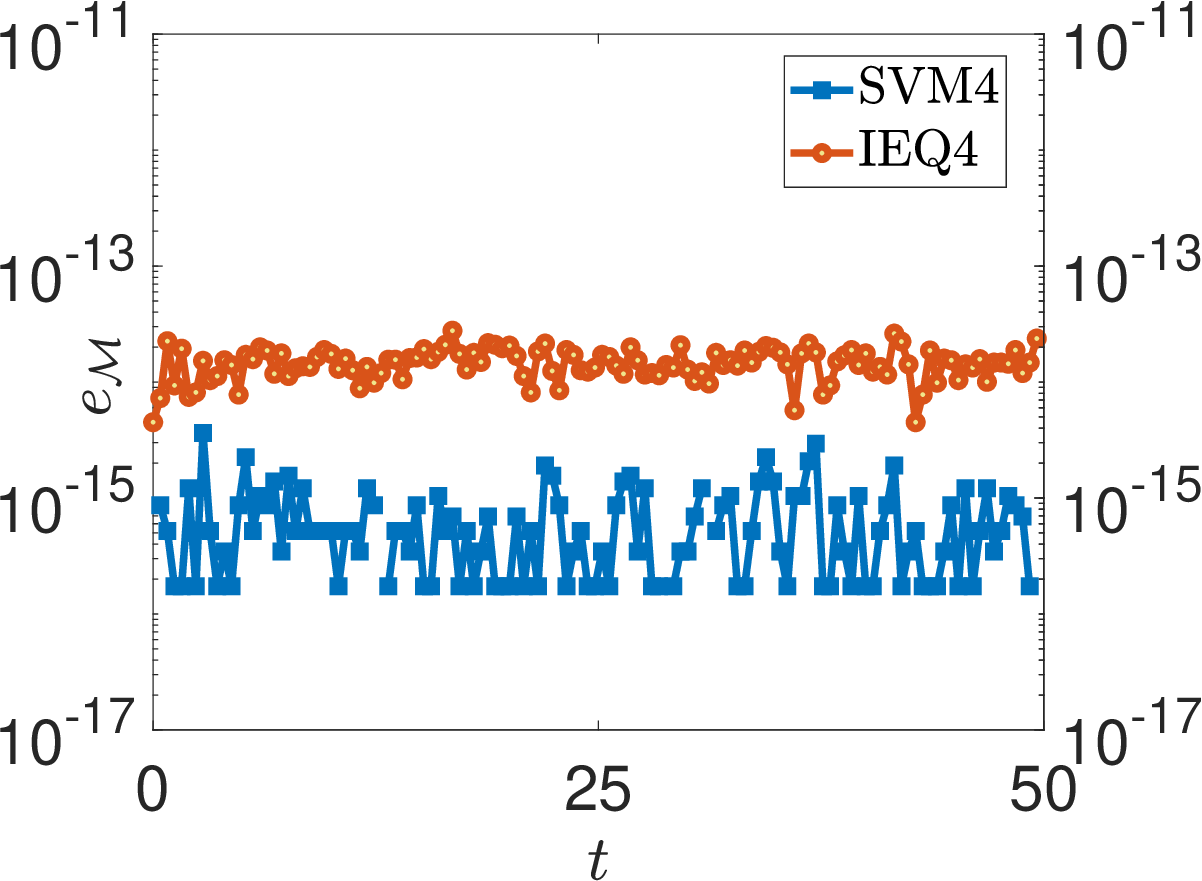}
\end{minipage}
\hspace{0.03\textwidth}
\begin{minipage}[t]{65mm}
\includegraphics[width=65mm]{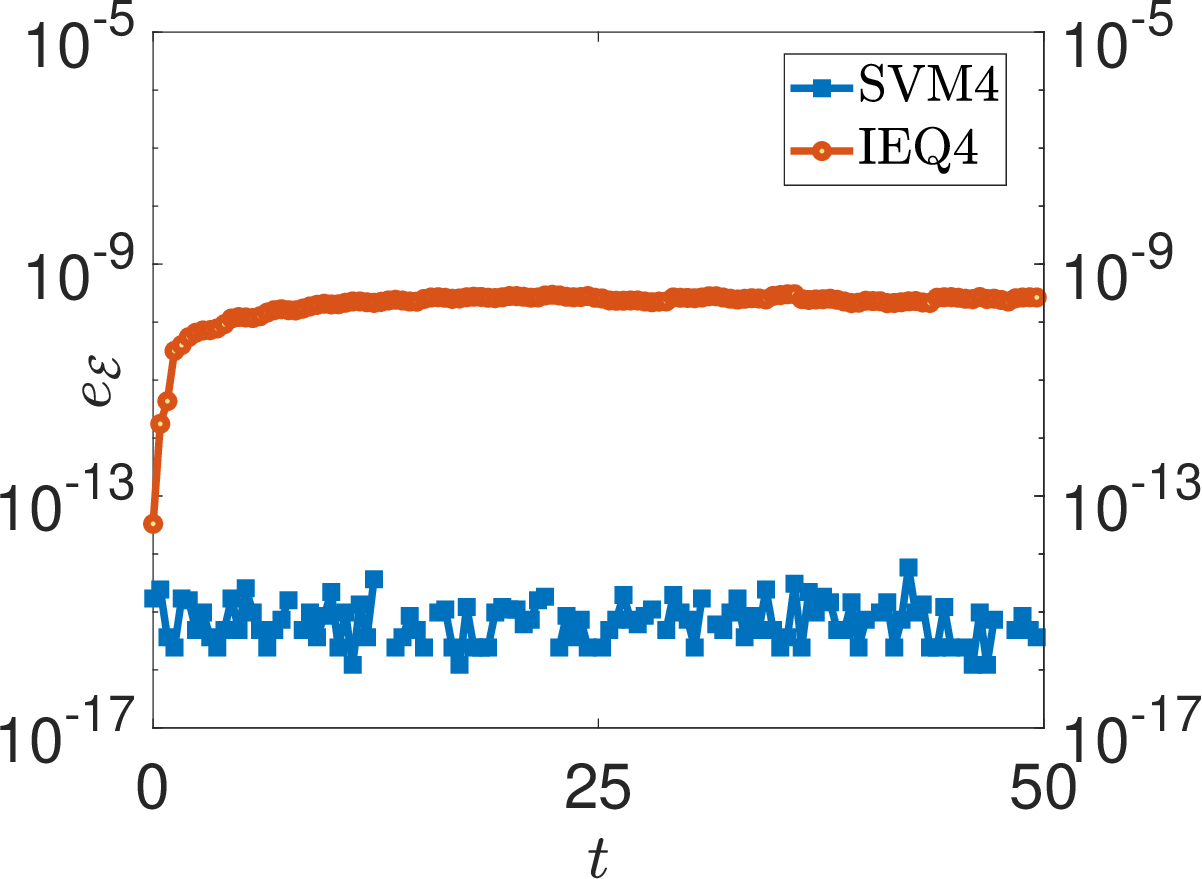}
\end{minipage}
\caption{Case II: The long-time evolution of $e_{\mathcal{M}}$ \& $e_{\mathcal{E}}$ for the two schemes in example \ref{EX4-4}.}\label{2d-scheme:case2}
\end{figure}

\begin{figure}[H]
\centering\begin{minipage}[t]{65mm}
\includegraphics[width=65mm]{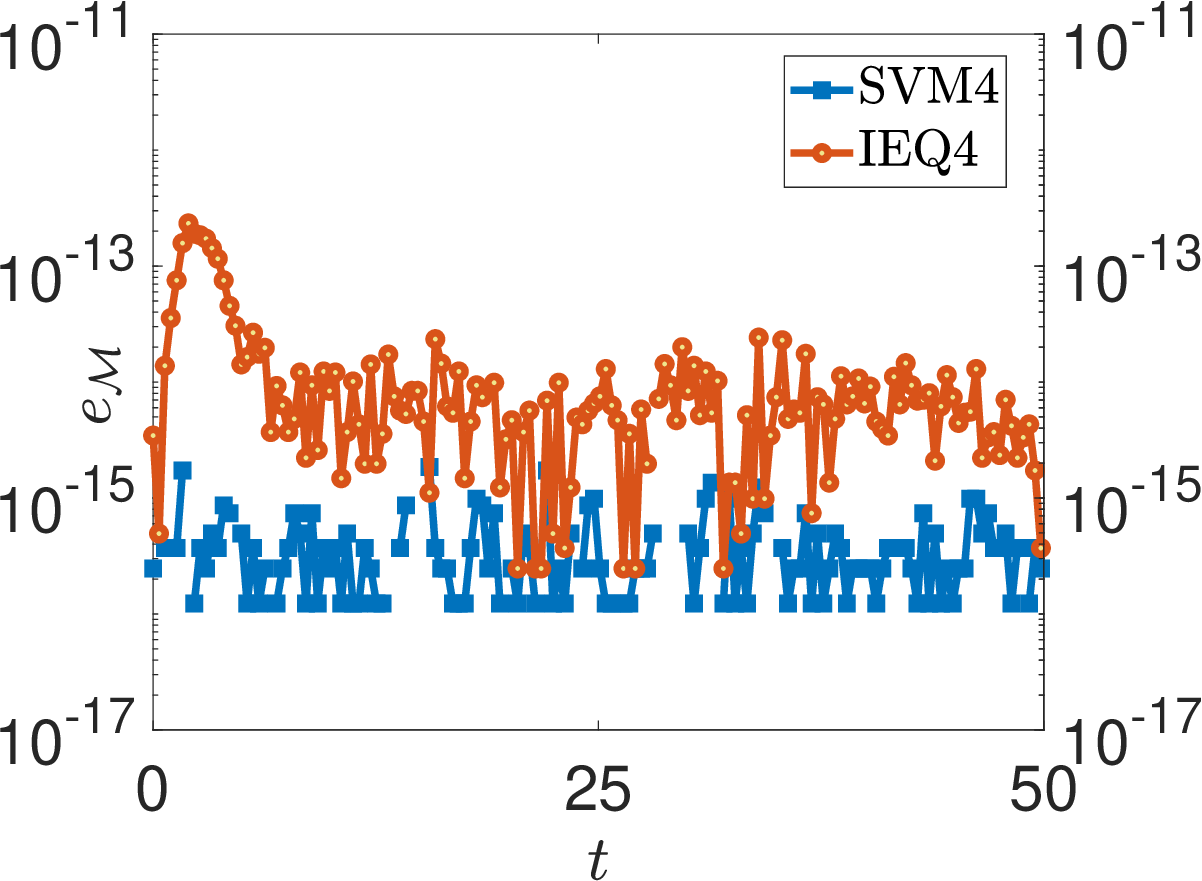}
\end{minipage}
\hspace{0.03\textwidth}
\begin{minipage}[t]{65mm}
\includegraphics[width=65mm]{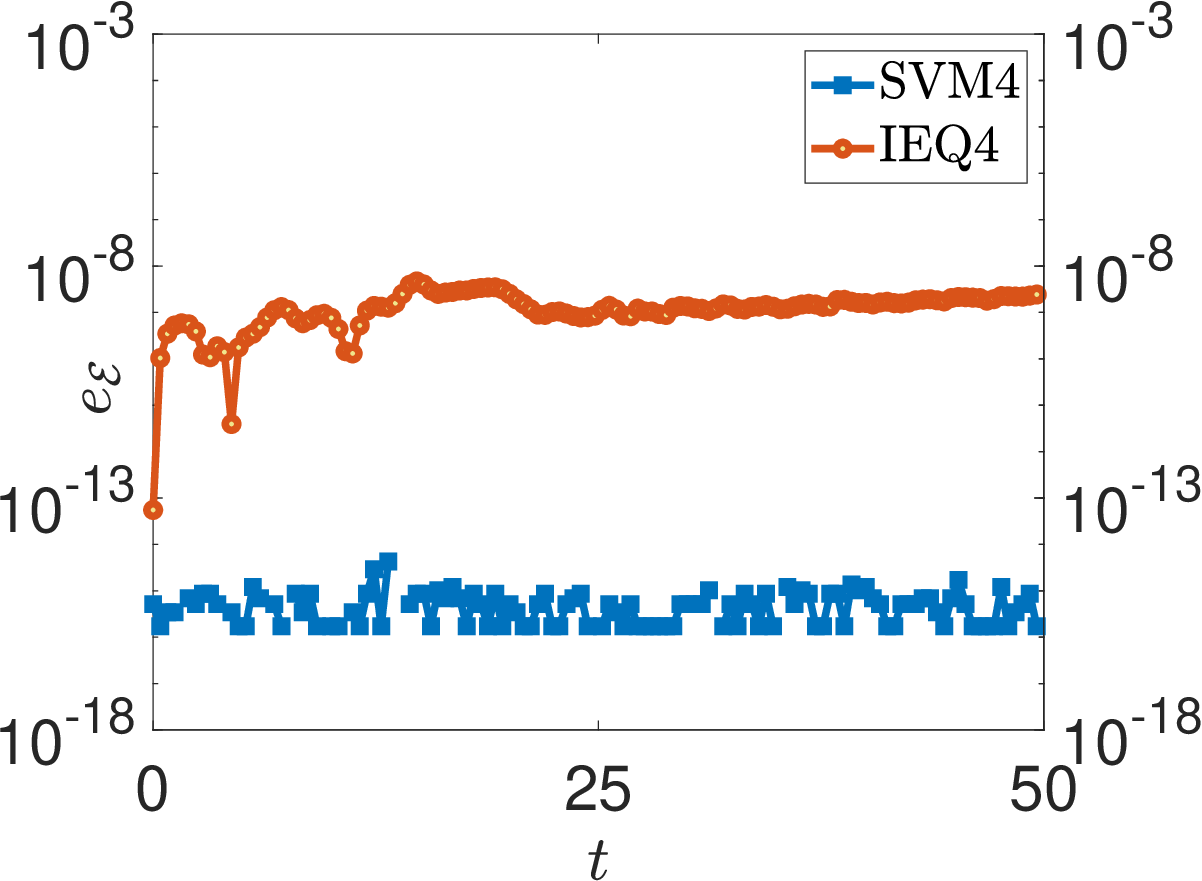}
\end{minipage}
\caption{Case III: The long-time evolution of $e_{\mathcal{M}}$ \& $e_{\mathcal{E}}$ for the two schemes in example \ref{EX4-4}.}\label{2d-scheme:case3}
\end{figure}

Figures \ref{2d-scheme:fig:11}-\ref{2d-scheme:fig:15} show the contour plots of $|u^{\epsilon}(\textbf{x},t)|^{2}$ at different time for Case I-III, respectively. Similar to the 1D case, we can observe from Figure \ref{2d-scheme:fig:11} that as two static Gaussons are in a sufficiently close proximity, they are attractive, collide and adhere momentarily before separating once again. It is worth noting that the Gaussons also oscillate in a pendulum-like motion, and small solitary waves are emitted outward during the interaction. Then, when velocities are introduced into one of them, Figure \ref{2d-scheme:fig:13} shows that the Gausson moving at a slower velocity will lead to the other moving in the same direction. However, if they are close enough, the moving Gausson will move perpendicular to the line connecting the two. While the static one is dragged to move, and the moving direction of the Gausson will be changed. Finally, two Gaussons will rotate and gradually move away from each other. This dynamics phenomena are consistent with the results obtained by Bao et al.~\cite{BCST-MMMAS-2022}.  Then, Figures \ref{2d-scheme:case1}-\ref{2d-scheme:case3} show the long-time evolution of $e_{\mathcal{M}}$ \& $e_{\mathcal{E}}$ for {\bf SVM4} and {\bf IEQ4} in example \ref{EX4-4}, from which we can see clearly that the errors on mass and energy produced by the {\bf SVM4} scheme can reach $10^{-15}$, while the {\bf IEQ4} scheme can only exactly preserve mass conservation law.

\section{Concluding remarks}\label{Sec:PM:5}
In this paper, we proposed a class of high-order, mass- and energy-preserving schemes for solving regularized logarithmic Schr\"{o}dinger equation. Based on the idea of the supplementary variable method, we first reformulate the original system into a new equivalent
system by introducing two supplementary variables, and then a fully-discrete scheme is presented by
using high-order prediction-correction scheme in time and Fourier pseudo-spectral method in space for
the reformulated system. The proposed schemes can preserve the original both mass and energy in discrete sense, and can reach arbitrary high-order accuracy in time. Additionally, in each time step, it only need to solve a constant-coefficient linear system accompanied by two scalar equations, which can be efficiently solved by the Newton iterative method. Numerical examples are addressed to illustrate the accuracy and conservation laws of the new method. Moreover, compared with the existing IEQ scheme, the newly developed method shows a remarkable advantage on mass and energy conservation laws for long-term numerical simulations. In the further work, we aim to establish optimal error estimates for the new method.
\section*{Acknowledgments}
This work is supported by the National
Natural Science Foundation of China
(Grant No. 12261097), the Yunnan Fundamental Research Project (Grant No. 202401AT070283), the Scientific Research Foundation Project of Yunnan University of Finance and Economics (Grant No. 2023C08) and the Graduate Student Innovation Foundation Project of Yunnan University of Finance and Economics (Grant No. 2024YUFEYC081).

\end{document}